\documentclass{amsart}
\usepackage{amssymb,amsmath,exscale,latexsym,amsthm,graphics,enumerate}  
\usepackage{epsfig}    % eps-Bilder einlesen
\usepackage{epic}      % irgendsone Graphikerweiterung
\usepackage{eepic}     % nochne erweiterung
\usepackage{subfigure} % mehre pics nebeneinander
\usepackage{psfrag}
\usepackage{booktabs}

\theoremstyle{plain}  \newtheorem{theorem}{Theorem}
                      \newtheorem*{theorem1_*}{Theorem 1'} 
                      \newtheorem*{theorem2_*}{Theorem 2'}
                      \newtheorem*{theorem*}{Theorem}
                      \newtheorem{lemma}{Lemma}[section]
                      \newtheorem{cor}[lemma]{Corollary}

\theoremstyle{definition}
        \newtheorem{definition}[lemma]{Definition}

\theoremstyle{remark}
        \newtheorem*{remark}{Remark}

        \newtheorem{problem}{Open Problem}

\theoremstyle{remark}

\theoremstyle{remark}

\theoremstyle{remark}
	\newtheorem*{claim}{Claim}

\numberwithin{equation}{section}
\newcommand{\dist}{\operatorname{dist}}

\newcommand{\const}{\operatorname{const}}

\newcommand{\diam}  {\operatorname{diam}}
\newcommand{\area}  {\operatorname{area}}

\newcommand{\inte}  {\operatorname{int}}

\newcommand{\id} {\operatorname{id}}
\newcommand{\R}{\mathbb{R}}      % reelle Zahlen
\newcommand{\C}{\mathbb{C}}      % komplexe Zahlen
\newcommand{\N}{\mathbb{N}}      % natürliche Zahlen
\newcommand{\Z}{\mathbb{Z}}      % ganze Zahlen
\newcommand{\Ho}{\mathbb{H}^+}   % obere Halbebene
\newcommand{\Hu}{\mathbb{H}^-}   % untere Halbebene
\newcommand{\Hobd}{\overline{\mathbb{H}}^+}   % abgeschlossene obere Halbebene
\newcommand{\Hubd}{\overline{\mathbb{H}}^-}   % abgeschlossene untere Halbebene
\newcommand{\CDach}{\widehat{\mathbb{C}}}% Riemannsche Zahlenkugel
\newcommand{\RDach}{\widehat{\mathbb{R}}}% aequator in CDach
\newcommand{\T}{\mathbb{T}}      % Torus
\providecommand{\defn}[1]{\emph{#1}}
\providecommand{\abs}[1]{\lvert#1\rvert}

\newcounter{mycount}
\newcommand{\SC} {\mathcal{S}} 

\newcommand{\OC} {\mathcal{O}}

\newcommand{\X} {\mathbf{X}}

\newcommand{\crit}{\operatorname{crit}}
\newcommand{\post}{\operatorname{post}}
\newcommand{\postpre}{\operatorname{post_p}}
\newcommand{\postper}{\operatorname{post_f}}

\begin{document}
 
\title{Dimension of elliptic harmonic measure of Snowspheres}

\author{Daniel Meyer}

\address{
    Department of Mathematics and Statistics\\
    P.O. Box 68 (Gustaf H\"{a}llstr\"{o}min katu 2b)\\
    FI-00014 University of Helsinki\\
    Finland}

\email{dmeyermail@gmail.com}

\thanks{This research was partially supported by an NSF postdoctoral
  fellowship, and NSF grant DMS-0244421.}

\keywords{Quasiconformal maps, elliptic harmonic measure, snowball,
  snowsphere, Latt\`{e}s maps,
  thermodynamic formalism} 

\subjclass[2000]{Primary: 30C65; Secondary 37F10, 37A05}

\date{\today}

\begin{abstract}
  A metric space $\mathcal{S}$ is called a \defn{quasisphere} if there is a
  quasisymmetric homeomorphism $f\colon S^2\to \mathcal{S}$. We consider the
  \defn{elliptic harmonic measure}, i.e., the push forward of
  $2$-dimensional Lebesgue measure by $f$. It is shown that for
  certain self similar quasispheres $\mathcal{S}$ (\defn{snowspheres}) the
  dimension of 
  the elliptic harmonic measure is strictly less than the Hausdorff
  dimension of $\mathcal{S}$. 
  This result is obtained by representing the self similarity 
  of a snowsphere by a
  postcritically finite rational map, and showing a corresponding
  result for such maps. As a corollary a metric characterization of
  Latt\`{e}s 
  maps is obtained. 
  Furthermore, a method to compute the dimension of elliptic harmonic
  measure 
  numerically is presented, along with the (numerically
  computed) values for certain examples. 
\end{abstract}

\maketitle

%\textpages

\section{Introduction}
\label{sec:introellip}
A homeomorphism $f\colon X\to Y$ of metric spaces is called
\defn{quasisymmetric} if 
there is a homeomorphism $\eta\colon [0,\infty)\to [0,\infty)$ such
that
\begin{equation*}
  \frac{\abs{f(x)-f(y)}}{\abs{f(x)-f(z)}}\leq \eta\left(\frac{\abs{x-y}}{\abs{x-z}}\right),
\end{equation*}
for all $x,y,z\in X, x\neq z$. Here and in the following, we use the
\defn{Polish notation} $\abs{x-y}$ for the metric. See \cite{JuhAn}
for general background on quasisymmetric maps. Every quasisymmetry is
trivially quasiconformal, see \cite{vaisala_n_dim} for background on
quasiconformal maps.

\medskip
A metric space $\SC$ is called a \defn{quasisphere} if it is
quasisymmetrically equivalent to the standard $2$-sphere $S^2$, i.e.,
if there is a quasisymmetry $f\colon S^2\to \SC$. To give a geometric
characterization of quasispheres is an open problem; the best results
to date are due to Mario Bonk and Bruce Kleiner \cite{BonKlei}. 

\medskip
In \cite{snowemb} and \cite{snowballquasiball} \defn{snowspheres} were
constructed, which are topologically $2$-dimensional analogs of
snowflake curves. These were shown to be quasispheres.
The construction is recalled briefly in the next section. Here, further
properties of qs-parametrizations of such snowspheres are investigated.
% The
% following is the main result of this paper. 

\begin{theorem}
  \label{thm:almostbiHoelder}
  Let $\SC$ be a snowsphere with uniformizing (quasisymmetric) map
  $f\colon S^2\to \SC$. There is $\alpha > 2/\dim_H(\SC)$ such
  that the following holds.
  For (Lebesgue) almost every $x\in S^2$, 
  \begin{equation*}
    \lim_{y\to x}\frac{\log\abs{f(x)-f(y)}}{\log\abs{x-y}}=\alpha.
  \end{equation*}
\end{theorem}
Note however, that the set where the above limit assumes a different
value/ does not exist is dense.  

\smallskip
Denote by $\lambda$ the $2$-dimensional, normalized ($\lambda(S^2)=1$)
Lebesgue measure on 
the sphere $S^2$. 
The \defn{elliptic harmonic measure} $\mu$ of $\SC$ induced by $f$ is
the pushforward of 
$\lambda$ by $f$,
\begin{equation}
  \label{eq:defmu}
  \mu=\mu_f=f_*\lambda.
\end{equation}
Recall that $f_*\lambda(A)=\lambda(f^{-1}A)$ for any Borel measurable
set $A\subset \SC$.

% We will restrict ourselves in this chapter to self similar snowballs/snowspheres, where the same generator is used throughout the construction. 
% Then we have a single rational function $R$ representing the self similarity as constructed in the last chapter.

% Here we consider the pullback of Lebesgue measure on the sphere under the uniformizing map $f\colon \SC\to S^2$.
% \begin{definition}
%   Given a quasiconformal map $f\colon \R^3\to \R^3$, the \emph{elliptic harmonic measure}\index{elliptic harmonic measure}\label{not:ellipmeas} of $D=f^{-1}(\B)$ is given by 
% $$\mu_f:=\lambda\circ f,$$
% where $\lambda$ denotes 
% Lebesgue measure on the sphere $S^2$, normalized by $\lambda(S^2)=1$.    
% \end{definition}

The \emph{dimension} of a probability measure $\mu$ is 
\begin{equation}
  \label{not:dimmu}
  \dim \mu:=\inf \{\dim_H(E) : \mu(E)=1\},
\end{equation}
where $\dim_H$ denotes the Hausdorff dimension.

\begin{theorem}
  \label{thm:main}
  The elliptic harmonic measure $\mu$ of every snowsphere $\SC$
  satisfies 
  \begin{equation*}
    \dim \mu =\frac{2}{\alpha} < \dim_H \SC.
  \end{equation*}
\end{theorem}

In concrete examples, the dimension $\dim \mu$ may be explicitly
computed numerically. We do this for several snowspheres.
 
The uniformizing map $f\colon S^2\to \SC$ is not unique for a given
quasisphere $\SC$. The dimension of the elliptic harmonic measure
however is independent of the particular map $f$.

\begin{lemma}
  \label{lem:dimellip}
  Let $f,g\colon S^2\to \SC$ be two quasisymmetric homeomorphisms. Then
  $$\dim \mu_f=\dim\mu_g.$$
\end{lemma}

\begin{proof} 
  The map 
  $$f^{-1}\circ g \colon S^2 \to S^2$$
  is quasisymmetric, and therefore maps zero sets (of Lebesgue measure
  $\lambda$ on the sphere $S^2$) to zero sets (see
  \cite{vaisala_n_dim}, Section 33). So $\mu_f$ and $\mu_g$
  are mutually absolutely continuous. 
\end{proof}
Hence, we can speak of \emph{the} dimension of elliptic harmonic
measure, which thus measures the size of the quasisphere
$\SC$ in a quasisymmetric fashion. A similar argument shows that
Theorem \ref{thm:almostbiHoelder} 
is independent of the specific quasisymmetric parametrization $f\colon
S^2\to \SC$. 

\medskip
Note that the analog statement for quasicircles is false. 
This is due to the fact that quasisymmetric maps $h\colon S^1\to S^1$
are not absolutely continuous in general \cite{AhlBeurext}. Thus, given two
quasisymmetric ho\-meo\-mor\-phisms $h_{1,2}\colon S^1\to \mathcal{C}$, one has in
general $\dim h_{1*}\lambda_1\neq \dim h_{2*}\lambda_1$ (here
$\lambda_1$ denotes $1$-dimensional Lebesgue measure on the circle $S^1$).
In fact, %it is not hard to show that 
the infimum of $\dim
h_{*}\lambda_1$ over all qs maps $h\colon S^1\to \mathcal{C}$ is
zero for all quasicircles $\mathcal{C}$ (see \cite{MR1051832}). Thus,
the notion of dimension 
of elliptic harmonic measure is unsuited to the $1$-dimensional
case. There are however interesting \emph{quantitative} questions in
this case, i.e., bounds on $\dim h_*\lambda_1$ depended on the
distortion of $h$.

\medskip
In \cite{MR2001c:49067}, surfaces $\SC\subset \R^3$ were constructed
which admit a parametrization $h\colon \R^2\to \SC$ satisfying
$1/C \abs{x-y}^\alpha\leq\abs{h(x)-h(y)}\leq C\abs{x-y}^\alpha$, for a
constant $C\geq 1$ and $0<\alpha<1$ (thus, $h$ is a quasisymmetry). It
follows that in these examples $\dim h_*\lambda =\dim_H \SC$. 

\medskip
In \cite{snowemb}, it was shown how the self similarity of a snowsphere
can be represented by a \defn{postcritically finite} rational map
$R$. The sphere $S^2$ can be 
equipped with a metric $\abs{x-y}_{\SC}$ with respect to which $R$
acts as a 
\defn{local similarity}. This means there is a $N>1$ such that for all
$x\in S^2$ there is a neighborhood $U(x)$ where
\begin{equation}
  \label{eq:def_local_sim}
  \frac{\abs{R(x)-R(y)}_\SC}{\abs{x-y}_\SC}= N,
\end{equation}
for all $y\in U(x)$ ($x\ne y$). In fact, $N$ is the ``scaling factor''
of the self 
similar snowsphere $\SC$. 
A piece of the snowsphere is isometric to a
hemisphere in this metric. 

Such metrics are constructed in \cite{BonkMeyer_MarkovThurston} for
all postcritically finite rational maps
that have no periodic critical points (equivalently, which have the
whole sphere as their Julia set). 
In fact, such metrics were constructed for \defn{expanding Thurston
  maps}, see below.
Theorem \ref{thm:almostbiHoelder} is
a consequence of the following.

\begin{theorem}
  \label{thm:Ralpha}
  Let $R$ be a postcritically finite rational map without periodic
  critical points. Let $\abs{x-y}_{\SC}$ be 
  a metric as
  above. Then there is $\alpha\geq 2/\dim_H(\SC)$ such that for
  (Lebesgue) almost every point $x\in S^2$ 
  \begin{equation*}
    \lim_{y\to x}\frac{\log \abs{x-y}_{\SC}}{\log\abs{x-y}} = \alpha.
  \end{equation*}
  Here, $\abs{x-y}$ denotes the spherical metric, $\dim_H(\SC)$ the
  Hausdorff dimension of $(S^2,\abs{x-y}_{\SC})$. Furthermore, $\alpha=
  2/\dim_H(\SC)$ if and only if $R$ is a Latt\`{e}s map.   
\end{theorem}
In fact, $\alpha=\frac{\log N}{\chi}$, where $\chi$ is the
\defn{Lyapunov exponent}. 

The next theorem (and hence Theorem \ref{thm:main}) is an easy
consequence. Let $\dim (\lambda,\SC)$ be the dimension of Lebesgue
measure with respect to the metric $\abs{x-y}_{\SC}$.  

\begin{theorem}
  \label{thm:Rmu}
  With the setup as above
  \begin{equation*}
    \dim(\lambda,\SC)=2/\alpha < \dim_H (S^2,\abs{x-y}_{\SC}),
  \end{equation*}
  unless $R$ is a Latt\`{e}s map, in which case there is 
  equality in the above.
\end{theorem}
Note that we can also express the dimension as
$\dim(\lambda,\SC)=\frac{2\chi}{\log N}=\frac{h}{\log N}$, where $h$
is the entropy.

\smallskip
In 
\cite{BonkMeyer_MarkovThurston}, \defn{expanding Thurston
  maps} were considered. These are maps that behave
\emph{topologically} as a rational map, 
meaning they can be written locally as $z\mapsto z^n$ after suitable
homeomorphic coordinate changes in domain and range. Furthermore, they
are \defn{postcritically finite} and satisfy an \defn{expansion}
property. Thurston has classified when such a map is \defn{equivalent}
to a rational map \cite{DouHubThurs}. 
For such an expanding Thurston map, a 
metric $\abs{x-y}_{\SC}$ satisfying (\ref{eq:def_local_sim}) was
constructed in \cite{BonkMeyer_MarkovThurston}, see also
\cite{HaiPilCoarse}. The metric is not 
unique, but two such metrics 
are snowflake equivalent. 
\begin{definition}
  Two metric spaces $(X_1,\abs{x-y}_1)$, $(X_2,\abs{x-y}_2)$ are
  \defn{snowflake equivalent} if there is a homeomorphism
  $\varphi\colon X_1\to X_2$, constants $C\geq 1, \beta>0$, such that
  \begin{equation*}
    \frac{1}{C}\abs{x-y}_1\leq \abs{\varphi(x)-\varphi(y)}^\beta_2\leq C\abs{x-y}_1,
  \end{equation*}
  for all $x,y\in X_1$.
\end{definition}

In
\cite{BonkMeyer_MarkovThurston} (see also \cite{HaiPilCoarse}), it is
proved that 
$(S^2,\abs{x-y}_{\SC})$ is quasisymmetrically equivalent to $S^2$
(with the standard metric) if and
only if $R$ is topologically conjugate to a rational map. In the present paper, we
obtain the following related theorem.

\begin{theorem}
  \label{thm:lipschitz_lattes}
  Let $R$ be an expanding Thurston map, and the metric
  $\abs{x-y}_{\SC}$ as in (\ref{eq:def_local_sim}).
  Then 
  \begin{align*}
    &(S^2,\abs{x-y}_{\SC}) \text{ is snowflake equivalent to
      (standard) } S^2
    \intertext{if and only if}
    &\text{$R$ is topologically conjugate to a
      Latt\`{e}s map.} 
  \end{align*}
\end{theorem}

We now give an outline of the paper.
In the next section, the construction of \defn{snowspheres} is recalled.

The construction of the rational map encoding the self similarity of a
snowsphere is done in 
Section \ref{sec:orig-with-rati}. This was only done for a specific
example in \cite{snowemb}, so full proofs are given. 

\smallskip
The proof of Theorem \ref{thm:Ralpha} becomes a relatively simple
application of the \defn{thermodynamic formalism}. The necessary facts
are reviewed in Section \ref{sec:ergdyn}.

\smallskip
The \defn{invariant measure} $\mu$ of $R$ that is absolutely
continuous with respect 
to Lebesgue measure is constructed in Section
\ref{sec:invar-meas-sphere}.  

\smallskip
The proofs of Theorem \ref{thm:Ralpha}, Theorem \ref{thm:Rmu}, and
Theorem \ref{thm:lipschitz_lattes} are
done in Section \ref{sec:estim-hold-expon}.  

\smallskip
In Section \ref{sec:jacob-invar-meas}, 
the strict inequality $\alpha > 2/\dim_H(\SC)$ is shown (unless $R$ is
a Latt\`{e}s map). This
already follows from \cite{ZdunikParabolicmaxmeas}, but can be shown
directly.

\smallskip
Section \ref{sec:numer-exper} shows how one can compute $\dim\mu$ using
Birkhoff's ergodic theorem. The results for several examples are
given.  

\smallskip
In the last section, some open problems are presented. 
% The harmonic measures of two complementary domains in $\R^3$ (disjoint with the same boundary) are 
% very different in general (singular). The elliptic harmonic measures for complementary domains are absolutely continuous by definition. Therefore one should expect harmonic measure and elliptic harmonic measure to be singular in general.  
% \begin{problem}
%   Is there a relation between the dimension of harmonic measure $\omega$ and the dimension of an elliptic harmonic measure $\mu$? Is $\dim \omega<\dim \mu?$
% \end{problem}

\subsection{Notation}
\label{sec:notation}

Two nonnegative quantities $A,B$ are \defn{comparable} if there is a
constant $C\geq 1$ such that
\begin{equation*}
  \frac{1}{C}A\leq B \leq C A.
\end{equation*}
We then write $A\asymp B$. The constant $C$ is referred to as
$C(\asymp)$. Similarly, we write $B\lesssim A$ if there is a constant
$C>0$ such that $B\leq C A$; we then refer to the constant $C$ by
$C(\lesssim)$. 

% The $\epsilon$-neighborhood of a set $S$ is denoted by
% $U(S,\epsilon)$.

The Riemann sphere is denoted by $\CDach=\C\cup\{\infty\}$. The
extended real line is $\widehat{\R}=\R\cup \{\infty\}\subset
\CDach$. The upper half plane is denoted by $\Ho$, the lower half
plane by $\Hu$; their closures by $\Hobd$, $\Hubd$. 

The $j$-th \defn{iterate} of a (rational) map $R$ is denoted by $R^j$.

\section{Snowspheres}
\label{sec:snowspheres}

Before presenting the general
construction, 
let us give an example first. This \emph{main} or \emph{standard}
example will serve throughout the whole paper to illustrate the
construction.  

Start with the unit cube. Divide each side into $5\times
5$ squares of side length $1/5$ ($1/5$\defn{-squares}). Put a cube of side
length $1/5$ (a $1/5$\defn{-cube}) on the middle $1/5$-square of each
side. This results in a polyhedron bounded by $6\times 29$ $1/5$-squares.
The construction is now iterated. Divide each $1/5$-square into
$5\times 5$ $1/25$-squares. Put a $1/25$-cube on the middle
$1/25$-square of each $1/5$-square, and so on. The limiting surface is called a
\defn{snowsphere} (which bounds a \defn{snowball}). 

In general, we divide each square into $N\times N$ $1/N$-squares
(where $N$ is an integer $\geq 2$). 
\begin{definition}
\label{def:generator}
  An  $N$\defn{-generator} $G$ is a polyhedral surface in $\R^3$ such that
  \begin{itemize}
  \item each face is a $1/N$-square.
  \item The boundary of $G$ is the unit square, $\partial G=\partial
    [0,1]^2\subset \R^3$. Here, we identify $\R^2\supset [0,1]^2$ with  
    $\R^2\times\{0\}\subset\R^3$.  
  \item $G$ is homeomorphic to $[0,1]^2$.
  \item The generator is \defn{symmetric} with respect to
    reflections on the planes $\{x=y\}, \{x+y=1\}, \{x=1/2\},
    \{y=1/2\}$, hence it is also symmetric with respect to rotations
    around the axis $\{x=y=1/2\}$ by multiples of $\pi/2$.
  \item Only one $1/N$-square in $G$ intersects each vertex of
    $[0,1]^2\subset \R^3$.  
  \end{itemize}
\end{definition}
This implies that all vertices of each $1/N$-square from which $G$ is
built are contained in the grid $\frac{1}{N}\mathbb{Z}^3$. 

The \defn{double pyramid} $\mathcal{P}$ is the union of the (solid)
pyramid with base $[0,1]^2\subset \R^3$ and tip
$(\frac{1}{2},\frac{1}{2},\frac{1}{2})$ and the (solid) pyramid with base
$[0,1]^2$ and tip $(\frac{1}{2}, \frac{1}{2}, -\frac{1}{2})$.  
 
A (self similar) snowsphere $\SC$ is constructed by repeatedly replacing
squares by scaled copies of the same generator. 
We exclude the trivial case where $G=[0,1]^2$.
The limiting surface
is a topological sphere embedded in $\R^3$ if the generator $G$ used in
the construction satisfies 
$G\cap \partial\mathcal{P} = \partial [0,1]^2$ (this implies that
$G\subset\mathcal{P}$, by the last property of $G$). 
It 
was shown in \cite{snowballquasiball} that it is quasisymmetrically
equivalent to the standard sphere $S^2$, here $\SC$ is equipped with
the metric inherited from $\R^3$.  

\subsection{Abstract Snowspheres}
\label{sec:abstract-snowspheres}

If the generator does not satisfy
$G\cap \partial\mathcal{P}=\partial[0,1]^2$, 
we can still
define the snowsphere generated by $G$ abstractly as follows. 
Choose the homeomorphism $h\colon G\to [0,1]^2$ from Definition
\ref{def:generator} such that it is the identity on $\partial
[0,1]^2$.  
View
$\SC_0:=\partial [0,1]^3$ as a two-dimensional cell complex in the
obvious way---faces, edges, vertices are the $2$-,$1$-, and
$0$-cells. In the same fashion, we view $G$ as a two-dimensional cell
complex. The cell complex $\SC_1$ is constructed by replacing each face
of $\SC_0$ by $G$. Formally, map the $2$-,$1$-, and $0$-cells of $G$ by
$h$ to $[0,1]^2$ and then to each face of $\SC_0$ by an
isometry. Images of cells of $G$ under this composition form the cell
complex $\SC_1$. 

This procedure is now iterated. Namely, map cells of $G$ by $h$ to
$[0,1]^2$ then to each $1/N$-square of $G$ by a similarity. 
So each $1/N$-square of the generator $G$ is a cell complex isomorphic
to $G$.
Map (the thus subdivided) $G$ to
$[0,1]^2$ by $h$ and subsequently to each side of $\partial [0,1]^3$
by an isometry  to obtain the cell complex $\SC_2$; and so on. Call 
the faces ($2$-cells) of $\SC_j$ the $j$\defn{-cylinders}. Note that
any $(j+1)$-cylinder $X_{j+1}$ is contained in exactly one
$j$-cylinder $X_j$. The $1$- and $0$-cells of $\SC_j$ are called
\emph{$j$-edges} and \emph{$j$-vertices}.

A more general treatment of \defn{subdivisions} can be found in
\cite{CFPfinSub}. It should be emphasized that the above 
procedure is used to define $\SC_j$ as a combinatorial object. The
metric on the abstract snowsphere will be defined exclusively from the 
combinatorics.

Consider
the set of sequences $x:=\{X_j\}_{j\geq 0}$, where each $X_{j+1}$ is a
$(j+1)$-cylinder contained in (the $j$-cylinder) $X_j$. Another such
sequence $y=\{Y_j\}$ is identified with $x$ if and only if
$X_j\cap Y_j\neq \emptyset$ for all $j$. The \defn{abstract
  snowsphere} $\SC$ generated by the generator $G$ is the set of all
equivalence classes of such sequences. 

The metric on $\SC$ is defined as follows. A $j$\defn{-chain} is a
sequence $Z_1,\dots,Z_n$ of $j$-cylinders such that $Z_k\cap
Z_{k+1}\neq \emptyset$, $k=1,\dots, n-1$. The \defn{length} of such a
$j$-chain is $n$. We say that the $j$-chain \defn{connects} the $j$-cylinders
$X_j$ and $Y_j$ if $Z_1=X_j$ and $Z_n=Y_j$. The chain is called
\defn{simple} if $Z_k\cap Z_m=\emptyset$ for $m > k+1$.          

Given two sequences $x=\{X_j\}, y=\{Y_j\}$ of $j$-cylinders as above
($X_{j+1}\subset X_j, Y_{j+1}\subset Y_j$), define
\begin{equation}
  \label{eq:defdn}
  d_j(x,y):=N^{-j}\min\{ \text{length of } j\text{-chain connecting }
  X_j, Y_j\}.
\end{equation}
Recall that $G$ was a $N$-generator, meaning it consists of
$1/N$-squares. 

\begin{theorem}
  \label{thm:djtod}
  For all $x,y$ as above, the limit
  \begin{equation*}
    d(x,y)=\lim_jd_j(x,y)
  \end{equation*}
  exists and is a metric on $\SC$.
\end{theorem}

\begin{proof}
  % test
  % \newcounter{mycounter}
  % \begin{list}{(\arabic{mycounter})}
  %   {\setlength{\leftmargin}{0pt}
  %     \setlength{\rightmargin}{0pt}
  %     \setlength{\itemsep}{0.8em}
  % %     \setlength{\labelwidth}{0pt}
  %     \usecounter{mycounter}}
  % \item  

  %   \label{item:dmetricfinite}
  (1)
  We first show that $d_j$ is \defn{finite}; more precisely $d_j(x,y)\lesssim
    1$, where $C(\lesssim)$ is independent of $x,y\in \SC$ and $j$.
    
    Consider the edges of length $1/N$ of the generator $G$; they are the
    edges of the $1/N$-squares from which $G$ is built. Note that each of
    the four edges of $\partial G=\partial [0,1]^2$ consists of
    precisely $N$ such edges.  
    
    This means that every $j$-edge can be covered by $N$
    $(j+1)$-edges. Hence by $N^k$ $(j+k)$-edges. 

    Let $M$ be the number of edges (of length $1/N$) in the
    generator, this is the number of $(j+1)$-edges in each
    $j$-cylinder $X_j$. Consider now a point $x=\{X_j\}\in \SC$. Then we
    can connect a given $0$-vertex $v_0$ of $X_0$ to a $1$-vertex $v_1$ of
    $X_1\subset X_0$ by at most $M$ $1$-edges. Similarly, we can connect
    $v_1$ to $X_2$ by at most $M$ $2$-edges and so on. Thus, $v_0$ can be
    connected to $X_j$ by at most 
    \begin{equation*}
      M(N^{j-1} + N^{j-2}+ \dots + 1) \lesssim N^j 
    \end{equation*}
    $j$-edges,
    where $C(\lesssim)=C(M,N)$. This shows the uniform boundedness of
    $d_j$.  

    \medskip
    For a given $j$-cylinder $X$, we define the \defn{annulus} (or
    $j$-annulus)  
    \begin{equation*}
      A(X):=\bigcup\{j\text{-cylinder }Y: Y\neq X, Y\cap X\neq
      \emptyset\}.  
    \end{equation*}
    Consider an $m$-chain $Z_1,\dots, Z_n$, where $m\geq j$. We say it
    \defn{crosses} 
    the annulus $A(X)$ if one $m$-cylinder $Z_{k}$ intersects
    $\partial X$, and another intersects the other boundary component of
    $A(X)$.
    \medskip

%  \item 
%    \label{item:dmetric_annulus}
    (2)
    A $(j+1)$-chain crossing a $j$-annulus has length at least $N$. 

    Consider first a chain of $1/N$-squares in the
    generator $G$ that connects two opposite (disjoint) sides of
    $\partial [0,1]^2$. This chain has length at least $N$. 
    
    The symmetry of the generator thus implies that a $(j+1)$-chain
    crossing a $j$-annulus has length at least $N$.
    \medskip
    %
    
%  \item 
%    \label{item:dmetric_conv}
    (3)
    We now 
    show that $d_j$ \emph{converges}.
    Let $x=\{X_j\}, y=\{Y_j\}\in \SC$ be arbitrary and distinct, $j_0\geq
    0$ be the smallest number such that
    $X_{j_0}\cap Y_{j_0}=\emptyset$. 
    Thus, 
    \begin{equation}
      \label{eq:dj0}
      d_{j_0}(x,y)\geq 3N^{-j_0}.     
    \end{equation}
    % Then any $j_0$-chain connecting
    % $x,y$ has to cross the annulus $A(X_{j_0}$, and $d_{j_0}(x,y) \geq
    % N^{-j_0}$. 
    We will show that $d_j(x,y) - 2 N^{-j}$ is increasing for $j\geq
    j_0$. This shows convergence of $d_j$ since it is bounded by (1). 
    % Let $x=\{X_j\}, y=\{Y_j\}\in \SC$ and $j_0$ be as above. 
    Let $j\geq j_0$. Consider a $(j+1)$-chain
    $Z_1,\dots,Z_n$ 
    connecting $X_{j+1}$ and $Y_{j+1}$. 
    % This chain has to cross the $j$-annulus $A(X_j)$, thus has length at
    % least $n\geq N$ by the above. Hence $d_{j+k}(x,y)\geq N^{-j}$.

    We
    construct a 
    $j$-chain $W_1,\dots, W_m$ connecting $X_j$ and $Y_j$ as follows. Let
    $W_1:=X_j$. Let $Z_{i_1}$ be the last $(j+1)$-cylinder that is
    contained in a $j$-cylinder $W_2$ intersecting $W_1$ (equivalently,
    the last $(j+1)$-cylinder in $A(W_1)$). The
    $(j+1)$-chain $Z_2,\dots,Z_{i_1}$ crosses the $j$-annulus
    $A(W_1)$. In the same fashion, let $Z_{i_2}$ be the last
    $(j+1)$-cylinder contained in a $j$-cylinder $W_3$ that intersects
    $W_2$. We continue to construct $\{W_i\}$ till $W_{m-1}$ intersects
    $W_m:=Y_j$. The chain $Z_2,\dots, Z_{n-1}$ crosses $m-2$ $j$-annuli
    $A(W_i)$. Thus, the length of the chain $\{Z_i\}$ is
    \begin{equation*}
      n\geq (m-2)N+2,
    \end{equation*}
    by (2). %(\ref{item:dmetric_annulus}). 
    Therefore
    \begin{align*}
      d_{j+1}(x,y) & \geq N^{-j-1}\left( (m-2)N +2\right)
      \\
      & \geq (m-2)N^{-j} + 2 N^{-j-1} \geq d_{j}(x,y) - 2N^{-j} + 2N^{-j-1}
      \\
      \Leftrightarrow \quad & d_j(x,y) -2N^{-j} \leq d_{j+1}(x,y) - 2N^{-j-1}. 
    \end{align*}
    Hence, $\lim_j (d_j(x,y) - 2 N^{-j}) = \lim_j d_j(x,y)=:d(x,y)$
    exists. 
    It does not degenerate since
    \begin{equation*}
      d_{j_0+k} (x,y) \geq d_{j_0}(x,y) - 2 N^{-j_0} + 2 N^{-j_0- k}\geq N^{-j_0},
    \end{equation*}
    for all $k\geq 0$ by (\ref{eq:dj0}). 
    The
    symmetry of $d$ is clear, the triangle inequality follows from the
    ones for $d_j$.   
    \medskip
    
    %
%  \item 
%    \label{item:dmetricwelldefined}
    (4)
    Finally, we show that the definition of $d(x,y)$ is \defn{independent
      of the representatives} $x,y$. 

    To verify this, let $x=\{X_j\} \sim
    \tilde{x}= \{\widetilde{X}_j\}$, $y=\{Y_j\}\sim
    \widetilde{y}=\{\widetilde{Y}_j\}$ (meaning that $X_j\cap
    \widetilde{X}_j\neq \emptyset$, $Y_j\cap \widetilde{Y}_j\neq
    \emptyset$ for all $j$). Consider a $j$-chain connecting
    $x$ and $y$. Adding $\widetilde{X}_j$ and $\widetilde{Y}_j$ to the
    beginning and end of this $j$-chain yields a $j$-chain connecting
    $\tilde{x}$ to $\tilde{y}$. Thus, 
    \begin{equation*}
      \abs{d_j(\tilde{x},\tilde{y})-d_j(x,y)}\leq 2 N^{-j}. 
      \qedhere
    \end{equation*}

%  \end{list}

%   Consider now a $(j+k)$-chain $Z_1, \dots , Z_n$ that
%   crosses the $j$-annulus 
%   $A(X)$. Assume that the chain is simple; furthermore that $Z_1$ is the only
%   cylinder intersecting $\partial X$, and $Z_n$ is the only cylinder
%   intersecting another boundary component of $A(X)$.
  
%   Let $Z_{i_1}$ be the last cylinder of the chain $\{Z_i\}$  that
%   is contained in 
%   a $(j+1)$-cylinder $Y_1$ which intersects $\partial X$. The
%   chain $Z_1,\dots Z_{i_1}$ crosses an annulus.

%   In the same fashion let $Z_{i_2}$ be the last cylinder of
%   $\{Z_i\}$ being contained in a $(j+1)$-cylinder $Y_2$ that
%   intersects $Y_1$. The chain $Z_{i_1},\dots, Z_{i_2}$ crosses the
%   $(j+1)$-annulus $A(Y_1)$. 
  
%   We end the procedure if the $(j+1)$-cylinder $Y_m$ intersects the
%   other boundary component of $A(X)$. Since $Y_1,\dots, Y_m$ crosses
%   $A(X)$, we have $m=N$. 

%   Thus $\{Z_i\}$ crosses at least $N$ $(j+1)$-annuli. The claim
%   follows by induction. 

\end{proof}
The proof above shows the following. If $m$ is the smallest number
such that $X_m\cap Y_m= \emptyset$ for given $x=\{X_j\}, y=\{Y_j\}$,
then $d(x,y)\asymp N^{-m}$. 

If $\SC$ is a snowsphere embedded in $\R^3$, then $d(x,y)$ is
comparable to the Euclidean distance of $x,y$ (see
\cite{snowballquasiball}, Section 2.4). 

The proof in \cite{snowballquasiball} applies here as well; each
abstract snowsphere is quasisymmetrically equivalent to the standard
sphere $S^2$.

\section{Origami with rational maps}
\label{sec:orig-with-rati}

In \cite{snowemb}, it was shown that the self similarity of a snowsphere
can be encoded by a rational map (this in turn was used to construct
the quasisymmetry to $S^2$). We review the construction briefly.

\subsection{The rational map representing the snowsphere.}
\label{sec:rati-map-repr}

\begin{figure}
  \centering
  \scalebox{0.9}{\includegraphics{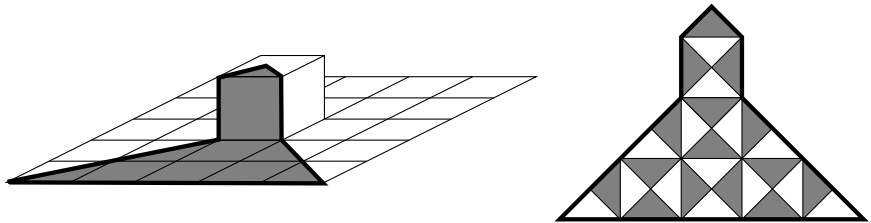}}
  \caption{Triangular generator of main example.}
  \label{fig:GT}
\end{figure}

Cut the generator along the diagonals $\{x=y\}, \{x+y=1\}$ into $4$
pieces. We call one such piece the \defn{triangular generator}
$G_T$. Figure \ref{fig:GT} shows the triangular generator of our main
example. 

\begin{lemma}
  \label{lem:triagGProp}
  The triangular generator $G_T$ satisfies the following.
  \begin{itemize}
  \item $G_T$ is compact, connected, and simply connected.
  \item Divide each $1/N$-square of $G$ into $4$
    $1/N$\defn{-triangles} along the diagonals. Then $G_T$ is built
    from such $1/N$-triangles. 
  \item Thus, $G_T$ may be viewed as a cell complex. The
    $1/N$-triangles, along with their edges/vertices are the
    $2$-,$1$-, and $0$-cells.
  \end{itemize}
\end{lemma}

\begin{proof}
  The first property follows from the symmetry of $G$. Recall that each
  $1/N$-square in $G$ has vertices in the grid
  $\frac{1}{N}\mathbb{Z}^3$. This implies the second property.   
\end{proof}
Note that each $1/N$-triangle has two vertices with angle $\pi/4$, and
one with angle $\pi/2$. 

There is a unique point of $G_T$ which intersects the line
$\{x=y=1/2\}$, which we call the \defn{tip} of the triangular
generator $G_T$. The two vertices of $\partial [0,1]^2$ contained in
$G_T$, together with the tip are the three \defn{corners} of $G_T$.
We think of $G_T$ as a topological
triangle. 

\begin{lemma}
  \label{lem:triagGProp2}
  Consider the vertices (of the $1/N$-triangles) in the triangular
  generator $G_T$.  
  \begin{itemize}
  \item Each of the two corners different from the tip is contained in
    exactly one $1/N$-triangle. The angle of the $1/N$-triangle at
    this corner is $\pi/4$. 
  \item If $N$ is odd the tip is contained in exactly one
    $1/N$-triangle, which forms an angle $\pi/2$ there. Otherwise, two
    $1/N$-triangles intersect at the tip, forming an angle $\pi/4$
    there. 
  \item The other vertices in $\partial G_T$ are contained in at least
    two $1/N$-triangles.
  \item Vertices not on $\partial G_T$ are contained in at least $4$
    $1/N$-triangles. 
  \end{itemize}
\end{lemma}

\begin{proof}
  The first property follows since only one $1/N$-square of the
  generator $G$ intersects each vertex of $[0,1]^2$. The other
  properties follow from the symmetry of $G$, the fact that $G$
  ``lives'' in the grid $\frac{1}{N}\Z^3$, and is homeomorphic to
  $[0,1]^2$.  
\end{proof}

Consider now two (identical) copies $G_T\times\{0\},G_T\times\{1\}$ of
the triangular generator.  
We now glue these two copies together along the boundary. More
precisely, we identify $x\times\{0\}, x\times\{1\}$ for 
each $x\in \partial G_T$. This yields a cell complex which is
topologically a $2$-sphere. 

We put a conformal structure on this
complex in the following way. Let $2\pi\alpha$ be the total angle at a
vertex $v$ (here each $2$-cell is viewed as a $1/N$-triangle). A chart
around $v$ is given by mapping all $2$-cells containing $v$ by
$z\mapsto z^{1/\alpha}$. 
The uniformization theorem gives a conformal
equivalence to the standard sphere $S^2$. We normalize the
uniformizing map 
by mapping the three corners of $G_T$ to
$-1,1,\infty$; where the tip is mapped to $\infty$. This implies 
that one copy of $G_T$ is mapped to the
upper half plane, while the other is mapped to the lower half
plane. Indeed, otherwise the normalized uniformizing map would fail to
be unique. So the common boundary of the two triangular generators
contained in our cell complex maps to the extended real line
$\RDach=\R\cup\{\infty\}$. 

For our standard example, the uniformizing map can be constructed
explicitly as follows. Map the polygon on the right in Figure
\ref{fig:GT} to the upper half plane by the inverse of a
Schwarz--Christoffel map. 

The images of the $2$-,$1$-,$0$-cells
are called $1$-\defn{tiles}, $1$-\defn{edges}, and
$1$-\defn{vertices} (these are all compact).  

\begin{lemma}
  \label{lem:tilesProp}
  The $1$-tiles have the following properties.
  \begin{itemize}
  \item The two $1$-tiles sharing a $1$-edge are conformal
    reflections of each other along this $1$-edge.
  \item At each $1$-vertex an even number of $1$-tiles intersect.
  \item Thus, we can color each $1$-tile either black or white, such
    that any two which share a $1$-edge are of different color.
  \end{itemize}
\end{lemma}

\begin{proof}
  The first property follows from the choice of charts when defining
  the conformal structure. 
 
  Consider a $1$-vertex $v\in \RDach$. Symmetry of the $1$-tiles with
  respect to $\RDach$ implies that an even number of $1$-tiles
  intersects in $v$. If $v\notin \RDach$, it is the image of a vertex
  $w\in G_T\setminus \partial G_T$. It is clear that at $w$ an even
  number of $1/N$-triangles intersect. This implies the second
  property. 

  The third property follows from the second.
\end{proof}

Color the $1$-tiles black and white as in the last lemma.
Recall that each $1/N$-triangle in the triangular generator $G_T$ had
exactly one vertex with angle $\pi/2$. Label all $1$-vertices which
are images of such vertices (of $1/N$-triangles) under the
uniformization by $b$. Label the two remaining $1$-vertices of each
white (black) $1$-tile $X'$ by $a,c$, such that $a,b,c$ are
mathematically 
positively (negatively) ordered in $\partial X'$. Note that the
labeling is compatible, meaning that each $1$-vertex gets exactly one
of the labels $a,b,c$ assigned.

Consider on a white $1$-tile $X'$ the Riemann map to the closed upper
half plane 
$\Hobd$, normalized by mapping $a\mapsto 1,b\mapsto \infty, c\mapsto
-1$. Let $Y'$ be a black $1$-tile sharing a $1$-edge with $X'$. By the
reflection principle,  the Riemann map $X'\to \Hobd$ extends
conformally to $Y'$, which
it maps to the closed lower half plane, again normalized by $a\mapsto
1,b\mapsto \infty, c\mapsto -1$.

Thus, we can define the \defn{rational map} $R$ \defn{representing the
  snowsphere} as follows. On each white (black)
$1$-tile $X'$ the map 
$R$ is the
Riemann map to the closed upper half plane $\Hobd$ (lower half plane $\Hubd$),
normalized by mapping 
$a\mapsto 1, b\mapsto \infty, c\mapsto -1$.
The map $R$ is well defined and holomorphic on $S^2=\CDach$ 
by the above (and the Riemann removability theorem). Hence, $R$ is a
rational map. 

There are two different (black/white)-colorings of $1$-tiles as
above---switch the colors of white $1$-tiles to black and vice
versa. The two rational maps $R,\widetilde{R}$ thus obtained satisfy
$\widetilde{R}= -R$. Both choices serve our intended purposes equally
well.  

The critical points are the $1$-vertices (where more than two
$1$-tiles intersect), they are
mapped by $R$ to either $1,\infty$, or $-1$ (if their label is
$a,b,c$).   
This allows to compute the map $R$ explicitly in many concrete examples. 
Namely, 
one obtains two expressions for $R$, in terms of $a$'s and $b$'s, as
well as in terms of $b$'s and $c$'s. 
Comparing
coefficients yields a system of equations which can be solved
numerically. For example, the rational map for our standard example
from Figure \ref{fig:GT} is
\begin{align*}
  & \widehat{R}(z)= 1+
  \\
  & \lambda 
  \frac{(z-1)(z-a_1)^4(z-a_2)^5(z+a_3)^3(z+a_4)^4(z-a_5)^4(z^2-a_6z+a_7)^4}
    {(z^2+t_1)^2(z^2+t_2)^2(z^2+t_3)^2(z^2-t_4)^2(z^2-t_5)^2(z^4-t_6z^2+t_7)^2}.
\end{align*}
Here,
\begin{align*}
  &\lambda=-0.001870\dots, 
  &t_1&=2712.82\dots,
  &t_2&=11.9805\dots,
  \\
  &t_3=0.31506\dots,
  &t_4&=1.01913\dots,
  &t_5&=3.97391\dots,
  \\
  &t_6=1.59735\dots,
  &t_7&=0.79867\dots,
  &a_1&=1.15921\dots,
  \\
  &a_2=9.00678\dots,
  &a_3&=191.820\dots,
  &a_4&=0.96246\dots,
  \\
  &a_5=0.48971\dots,
  &a_6&=-1.98760\dots,
  &a_7&=1.99159\dots \quad.
\end{align*}
In \cite{snowemb}, several other explicit examples of 
rational maps encoding certain snowspheres can be found.

Denote the set of critical points of $R$ by $\crit(R)$. 
%Recall that the set of critical points of $R$ are the $1$-vertices,
%$\crit (R)=\{1\text{-vertices}\}$. They are
%mapped to $\{1,\infty,-1\}$, which are again $1$-vertices. 
The \defn{postcritical set} is
\begin{equation*}
  \post(R):= \bigcup_{j\geq 1}R^j(\crit(R)).
\end{equation*}
%is $\post(R)=\{1,\infty,-1\}$.

%In particular $R$ is \defn{postcritically finite}, meaning that
%$\post(R)$ is a finite set.

\begin{lemma}
  \label{lem:propR}
  The map $R$ is \defn{postcritically finite}, meaning that $\post(R)$
  is a finite set. In fact $\post(R)=\{1,\infty,-1\}$. 
  Furthermore, $R$ has \emph{no critical periodic orbit},
  meaning that 
  $$R^k(z_0)=z_0\Rightarrow z_0\notin \crit(R),$$
  for $k\geq 1$.  
\end{lemma}

\begin{proof}
  Recall that the number of $j$-tiles intersecting in a $j$-vertex is
  even (Lemma \ref{lem:tilesProp}). 
  \begin{claim}
    A point $c$ is a critical point if and only if $c$ is a $1$-vertex
    at which at least four $1$-tiles intersect. 
  \end{claim}
  Let $c$ be a $1$-vertex as above, then $R$ maps tiles containing $c$
  alternatingly to the upper and lower half plane. Thus, $c$ is a
  critical point. If only two $1$-tiles intersect in the $1$-vertex
  $v$ it is 
  not a critical point by the same argument.

  To see the other implication, let $x\in \CDach$ be not a
  $1$-vertex. Then $x$ is contained in a set of the form $\inte X'\cup
  \inte Y'\cup \inte E$. Here, $X',Y'$ are $1$-tiles sharing the
  $1$-edge $E$. The map $R$ maps this set conformally (to either
  $\Ho\cup\Hu \cup (-1,1), \Ho\cup \Hu \cup (1,\infty)$, or $\Ho\cup
  \Hu\cup (-\infty,-1)$). So $x$ is not a critical point. 

  \medskip
  The map $R$ was constructed to map all $1$-vertices (hence, all
  critical points) to $1,\infty$, or $-1$. Since these are in turn
  $1$-vertices the map is postcritically finite. 

  \medskip
  It remains to show that $R$ has no critical periodic orbits. Recall
  that 
  $1,-1$ are the images of the two corners of the triangular generator
  $G_T$ different from the tip. By Lemma \ref{lem:triagGProp2}, two
  $1$-tiles intersect at both $1,-1$. Thus, they are not critical
  points. Again, by Lemma \ref{lem:triagGProp2}, the label at $1,-1$ is
  $a$ or $c$. Thus, $R$ maps $\{1,-1\}$ to $\{1,-1\}$.

  \medskip
  Consider now $\infty$, which is the image of the tip of the
  triangular generator $G_T$ under the uniformization. If $N$ is odd, 
  two $1$-tiles contain $\infty$, and $\infty$ is labeled $b$ (Lemma
  \ref{lem:triagGProp2}). Thus, $\infty$ is not a critical point and
  $R(\infty)=\infty$.

  If $N$ is even, four $1$-tiles contain $\infty$, and $\infty$ is
  labeled $a$, or $c$. Thus, $\infty$ is a critical point and mapped to
  $1$, or $-1$ by $R$.

\end{proof}

\subsection{Embedding the Snowsphere via the rational map}
\label{sec:embedd-snowsph-via}
We briefly describe in what sense the rational map $R$ ``encodes'' the
self similarity of the snowsphere. We equip the sphere $S^2$ with a
metric such that the upper and lower half planes $\Ho,\Hu$ are
each isometric to a ``triangular piece'' of the snowball. 

\smallskip
Consider one face ($0$-cylinder) $X$ of the snowball $\SC$. Let
$E\subset X$ be a $0$-edge, $E'=\partial X\setminus E$ be the union of
the other $0$-edges of $X$. A \defn{triangular piece} of the snowball
is defined as
\begin{equation*}
  \SC_T:=\{ x\in X \mid \dist(x,E) \leq \dist(x, E')\}. 
\end{equation*}

A $j$\defn{-tile} is one component of the preimage of a $1$-tile by
$R^{j-1}$ (equivalently one component of the preimage of $\Hobd$,
or $\Hubd$ by $R^j$). Preimages of $\{1,\infty,-1\}$ by $R^j$ are
called $j$-vertices. Note that each $(j+1)$-tile is contained in
exactly one
$j$-tile, since each $1$-tile is contained in either $\Hobd$ or
$\Hubd$. 

\begin{lemma}
  \label{lem:jtilesshrink}
  The size of $j$-tiles goes to zero,
  \begin{equation*}
    \max \diam X'_j \to 0 \quad \text{as } j\to \infty.
  \end{equation*}
  Here, the maximum is taken over all $j$-tiles $X'_j$. 
\end{lemma}

\begin{proof}
  This is a straightforward application of Schwarz' lemma. We refer the
  reader to \cite{snowballquasiball} Lemma 3.3 for details. 
\end{proof}
We could use the $j$-tiles to define a metric
on $S^2$ similarly as in Theorem \ref{thm:djtod}. 

We choose a slightly different approach here. Label all
$j$-vertices that are preimages of $\infty$ by $b$. Form the union of
$j$-tiles that intersect in a $j$-vertex labeled $b$ to form a
$j$\defn{-quadrilateral}. 

If $N$ is odd there is a single $j$-quadrilateral containing only two
$j$-tiles, namely the one containing $\infty$ (see Lemma
\ref{lem:triagGProp2} and the proof of Lemma \ref{lem:propR}). With
this exception, all $j$-quadrilaterals contain exactly $4$ $j$-tiles. 

Use the $j$-quadrilaterals to define a
metric as in Section \ref{sec:abstract-snowspheres},
\begin{align*}
  d'_j(x,y):= & N^{-j}\min\{\text{length of chain of }
  \\
  & j\text{-quadrilaterals connecting }x,y\}.   
\end{align*}

\begin{theorem}
  \label{thm:dS2}
  The limit $d'_j\to d'$ exists and is a metric on $S^2$. Furthermore,
  \begin{enumerate}
  \item \label{item:dS2isometric}
    $(\SC_T,d)$ is isometric to $(\Hobd, d')$ (as well as $(\Hubd, d')$). 
  \item \label{item:dS2qs}
    The map $(S^2,\abs{x-y}) \to (S^2,d')$ is
    \defn{quasisymmetric}. Here, 
    $\abs{x-y}$ is the spherical metric on $S^2$. 
  \item \label{item:dS2Rsimil}
    The map $R\colon (S^2,d')\to (S^2,d')$ is a \defn{local
      similarity}. More precisely, for any $x\in S^2$ there is a
    neighborhood $U$ of $x$ such that
    \begin{equation*}
      \frac{d'(R(x),R(y))}{d'(x,y)}=N,
    \end{equation*}
    for all $y\in U$.
  \item \label{item:HdimS2} 
    The \defn{Hausdorff dimension} of $(S^2,d')$ is
    \begin{equation*}  \label{eq:1dimStilde}
      \dim_H(S^2,d')=\frac{\log \deg R}{\log N}.
    \end{equation*}
%     Note that $\deg R$ is the number of $1/N$-squares in the generator
%     $G$. 
  \end{enumerate}

\end{theorem}

\begin{proof}
  (\ref{item:dS2isometric}) 
  Observe first that for $z,w\in \Ho$ a shortest connecting chain
  of 
  $j$-quadrilaterals is always in $\Ho$ (for $j$ sufficiently large).
  The same is true for $x,y\in \inte \SC_T$. 
  Note that the combinatorics of $j$-quadrilaterals in $\Ho$ is the
  same as the 
  one for the $j$-cylinders in $\inte \SC_T$. Map points
  $x,y\in \inte 
  \SC_T$ to points
  $z,w\in \Ho$ with the same combinatorics, then
  $d_j(x,y)=d'_j(z,w)$. The case when one or both points lie on the
  boundary is left to the reader. 
   
%   Let $\widetilde{X}_1\subset S^2$ be a
%   $1$-quadrilateral as in the statement, $X_1\subset\SC$ be a
%   $1$-cylinder.  
%   Note that the combinatorics of $j$-quadrilaterals contained
%   in $\widetilde{X}_1$ is the same as the combinatorics of
%   $j$-cylinder contained in $X_1$ by our
%   construction. Thus we can map $X_1$ to $\widetilde{X}_1$ by
%   mapping each $j$-cylinder $X_j$ to the $j$-quadrilateral with the same
%   combinatorics.
%   Note that for two points $x,y\in X_1$ (or in $\widetilde{X}_1$), 
%   (one of) the shortest $j$-chain connecting $x,y$ is
%   contained in $X_1$ (in $\widetilde{X}_1$).
 
%   Thus $d_j$ converges to a metric $d$ on $S^2$, by
%   Theorem \ref{thm:djtod}. The map
%   $X_1\subset(\SC,d) \to \widetilde{X}_1\subset (S^2,d')$ is an isometry.      
  
  \medskip
  From (\ref{item:dS2isometric}),
  Lemma \ref{lem:jtilesshrink}, and the
  proof of Theorem \ref{thm:djtod} it follows that $d'$ is a metric on
  $S^2$. 

  \medskip
  The proof for (\ref{item:dS2qs}) is identical to the one in
  \cite{snowemb}, and will not be repeated here. It follows however
  essentially from the facts recorded in the next lemma.

  \medskip
  (\ref{item:dS2Rsimil}) Let $x\in S^2$. Assume first that $x\neq
  \infty$, and that $x$ is not a pole of $R$. The union of
  $j$-quadrilaterals containing $x$ forms a neighborhood $U$ of
  $x$. Assume $j$ is large enough, such that $\infty\notin U$ and $U$
  does not contain a pole of $R$. Then every $j$-quadrilateral
  $\widetilde{X}_j$ (and every $(j+k)$-chain in 
  $\widetilde{X}_j$) is mapped by $R$ bijectively
  to a $(j-1)$-quadrilateral $\widetilde{X}_{j-1}=R(\widetilde{X}_j)$
  (to a $(j+k-1)$-chain).  
  As before, any shortest $(j+k)$-chain connecting $x,y\in
  \widetilde{X}_j$ can be chosen to lie in $\widetilde{X}_j$;
  similarly any shortest $(j+k-1)$-chain connecting $R(x),R(y)\in
  \widetilde{X}_{j-1}$ can be chosen to lie in
  $\widetilde{X}_{j-1}$.  
  The statement follows.

  The case where $x$ is either $\infty$ or a pole of $R$ are left to
  the reader. 

  \medskip
  (\ref{item:HdimS2})
  Since $(S^2,d')$ satisfies the \emph{open set condition}, its
  Hausdorff dimension is given by the \defn{pressure formula} (see
  \cite{Falconer1}, 
  p. 118, Theorem 9.3).  

\end{proof}

Property (\ref{item:dS2Rsimil}) shows how $R$ encodes the self
similarity of the snowsphere (each $j$-cylinder $X_j\subset \SC$ is
similar to a $(j+k)$-cylinder, scaled by the factor $N^{-k}$).  

It is easy to construct the quasisymmetric uniformizing map
$f\colon S^2\to \SC$ from the above. 
Namely uniformize the (surface of the) cube as in Section
\ref{sec:rati-map-repr}. Map a $1$-quadrilateral conformally into each
of the uniformized faces. The combinatorics of the
$(j+1)$-quadrilaterals is identical to the combinatorics of the
$j$-cylinders, so each point on $S^2$ can be mapped to a point on
$\SC$. The last theorem shows that this map is a quasisymmetry. 

\begin{cor}
  \label{cor:3implies1_4implies2}
  Theorem \ref{thm:Ralpha} implies Theorem \ref{thm:almostbiHoelder},
  and Theorem \ref{thm:Rmu} implies Theorem \ref{thm:main}. 
\end{cor}

From now on, we will denote the metric $d'$ from the last theorem by
$\abs{x-y}_\SC$. The spherical metric is denoted by $\abs{x-y}$.

\medskip
Using the $j$-tiles, it is possible to measure distances in purely
combinatorial terms. 
Let
$x,y\in S^2$ be arbitrary. Then 
\begin{equation}
  \label{eq:defjxy}
  j(x,y):= \min \{j \mid \text{ there exist disjoint }
  j\text{-tiles } X_j\ni x, Y_j\ni y\}.
\end{equation}
This measures to what level of tiles one needs to descend to be able to
distinguish $x$ and $y$. The following are from
\cite{snowballquasiball} (Lemma 3.10, Lemma 2.4, Lemma 3.7, and
Corollary 3.8).  
\begin{lemma}
  \label{lem:combmetric}
  For $x,y\in S^2$, $j=j(x,y)$ as above
  \begin{align}
    \label{eq:dS2xy}
    &\abs{x-y}\asymp \diam X'_j\\
    \label{eq:dSxy}
    &\abs{x-y}_{\SC}\asymp N^{-j},
  \end{align}
  where $C(\asymp)=C(N)$, and $X'_j$ is a $j$-tile containing
  $x$. 
  Furthermore, (for $j$-tiles $X'_j,Y'_j$)
  \begin{align}
    \label{eq:diamXY}
    &\diam X'_j\asymp \diam Y'_j, \quad \text{ if } X'_j\cap Y'_j\neq
    \emptyset, 
    \\
    \label{eq:diamXX}
    &\diam X'_{j+1} \asymp \diam X'_j, \quad \text{ for any $(j+1)$-tile
    } X'_{j+1}\cap X'_j\neq \emptyset,
    \\
    \label{eq:diamXareaX}
    &\area(X'_j)\asymp (\diam X'_j)^2,
  \end{align}
  with a constant $C(\asymp)$. 
\end{lemma}
Here, diameter and area are measured with respect to the spherical
metric. 

\subsection{Expanding postcritically finite maps.}
\label{sec:expand-postcr-finite}
Let $R$ be a postcritically finite rational map that has no critical
periodic orbit.
In \cite{BonkMeyer_MarkovThurston}, a metric $\abs{x-y}_\SC$ with
respect to which $R$ 
is a local similarity (as in Theorem \ref{thm:dS2}
(\ref{item:dS2Rsimil})) was
constructed. See also \cite{HaiPilCoarse} and \cite{CFKP}.
The map $\id \colon(S^2,\abs{x-y})\to (S^2,\abs{x-y}_\SC)$ is again
quasisymmetric. 
The Hausdorff dimension of $(S^2,\abs{x-y}_{\SC})$ is again given by
the formula in Theorem \ref{thm:dS2} (\ref{item:HdimS2}).
The metric $\abs{x-y}_\SC$ (as well as the expansion
factor $N$) is not unique, two such metrics are snowflake equivalent. 

\smallskip
Distances may again be measured in purely combinatorial terms. Pick a
Jordan curve $\mathcal{C}$ through $\post(R)$, we require that
$\mathcal{C}$ 
is a zero set with respect to 
$2$-dimensional Lebesgue
measure. The closure of one 
component of $R^{-j} (S^2\setminus \mathcal{C})$ is called a
$j$-tile. Lemma \ref{lem:combmetric} continues to hold. By
\cite{BonkMeyer_MarkovThurston} we can choose the curve
$\mathcal{C}\supset \post(R)$ to be invariant for some iterate $R^n$,
meaning that $R^n(\mathcal{C})\subset \mathcal{C}$. This means that
every $j$-tile defined in terms of $R^n$ and $\mathcal{C}$ is
contained in exactly one $(j-1)$-tile.

\subsection{The Orbifold associated to $R$.}
\label{sec:orbif-assoc-r}

We remind the reader of the notion of the \defn{orbifold}
  associated to a postcritically finite rational map $R$. See
\cite{MR2193309} (Appendix E) and \cite{McM}
(Appendix A) for more 
background. The \defn{ramification function} $\nu\colon \CDach\to
\N\cup\{\infty\}$ is 
defined by the following. It is the smallest function such that
\begin{align*}
  & \nu(z)=1, \text{ if } z \in \CDach\setminus \post(R),
  \intertext{and for all $p\in \post(R),R(q)=p$}   
  & \nu(p) \text{ is a multiple of } \deg_R(q)\nu(q).
\end{align*}
Clearly, $\nu$ is finite on $S^2$ if and only if $R$ has no critical
periodic orbit (see Lemma \ref{lem:propR}).  
We call $\OC=(S^2,\nu)$ the \defn{orbifold} associated to $R$. 
The \defn{signature} of $\OC$ is the list of values of $\nu$ at the
postcritical points. 
The
\defn{Euler characteristic} of $\OC$ is given by
\begin{equation*}
  \chi(\OC) = 2 + \sum_{p\in\post(R)} \left(\frac{1}{\nu(p)} - 1
    \right).  
\end{equation*}
It is nonpositive. The orbifold $\OC$ is called \defn{parabolic} if
$\chi(\OC)=0$, hyperbolic otherwise. 
Parabolic orbifolds are completely classified
(\cite{DouHubThurs}, Section 9). Indeed one checks directly that the
only possible signatures are $(\infty,\infty)$, $(2,2,\infty)$,
$(2,2,2,2)$, $(2,4,4)$, $(2,3,6)$, $(3,3,3)$. The last four of those
belong to rational maps $R$ having no critical periodic orbit. Such a
map is a \emph{Latt\`{e}s map}, i.e., obtained as a quotient of a map
from a torus $\T^2=\C/\Lambda$ to itself, $\varphi \colon \T^2 \to
\T^2$ (see \cite{MR2348953}, \cite{MR2193309}\footnote{Note that in
  earlier editions the definition 
  of Latt\`{e}s maps differs.}, and \cite{Lattes}).

Consider the ``trivial'' snowsphere, i.e., the one with
generator $G=[0,1]^2$. One checks directly that in this case the
corresponding rational map $R$  
is a Latt\`{e}s example with signature $(2,4,4)$. 

In all other cases, there are vertices in the generator where not $4$
$1/N$-squares intersect. The signature of the corresponding rational
map will contain at least one number $\geq 12$, thus $R$ is not a
Latt\`{e}s map.

\section{Ergodic Theory and Dynamics}
\label{sec:ergdyn}
General background for ergodic theory can be found in \cite{Walters}
and \cite{Petersen}. A survey of the methods used is given in
\cite{urbanski}, the forthcoming book \cite{przyurban} will contain an
exhaustive treatment.   The booklet by Michel Zinsmeister
\cite{Zinsmeister} is very readable, though (or possibly because) it
only deals with the hyperbolic case, not the subhyperbolic one at
hand.

Let $R\colon \CDach \to\CDach$ be a postcritically finite rational
map with no critical periodic orbit.
The reader should think of $R$ as the one
 constructed in  
Section \ref{sec:orig-with-rati}, which represents the self similarity
of a snowsphere and embeds it. 
In this case the postcritical set is
$\post(R)=\{-1,1,\infty\}$.
% ; some of the examples from \cite{snowemb}
% have different (though finite) postcritical sets. The map $R$ has no
% critical periodic orbits (see Lemma \ref{lem:propR}).

We consider two metrics on the sphere $S^2$. The ``self similar''
metric $\abs{x-y}_\SC$ is the one from Theorem 
\ref{thm:dS2} or from Section \ref{sec:expand-postcr-finite}.
The standard 
spherical metric is denoted by $\abs{x-y}$. When writing 
$\CDach$, we always mean the sphere equipped with the spherical metric. 
Theorem
\ref{thm:almostbiHoelder} will follow from Theorem \ref{thm:Ralpha} by
Theorem \ref{thm:dS2} (\ref{item:dS2isometric}). 
% Instead of the quasisymmetric uniformizing map $f\colon S^2\to \SC$,
% we will consider the map $F\colon \CDach=(S^2,\abs{x-y}) \to
% (S^2,\abs{x-y}_\SC)$, where 
% $F=\id_{S^2}$. This is enough 
% by Section \ref{sec:embedd-snowsph-via}, meaning that $\dim
% F_*\lambda_2=\dim f_* \lambda_2$. 

\smallskip
We denote 
the set of all $j$-tiles by $\X'_j$. Recall that if $R$ is a rational
map representing a snowsphere (Section \ref{sec:rati-map-repr})
they are given as the set of preimages of $\Hobd$ and $\Hubd$ under
the iterated map $R^j$.  
 Let $\mu$ be an $R$-\emph{invariant}\index{invariant
  measure}\label{not:muprime} probability measure, meaning that 
\begin{equation}\label{eq:invmeas}
  R_*\mu(A)=\mu(R^{-1}A)=\mu(A)
\end{equation}
for all Borel sets $A\subset\CDach$. Equivalently, for all $g\in L^1(\mu)$
\begin{equation}\label{eq:invmeas2}
  \int g\circ R\,d\mu=\int g\,d\mu.
\end{equation}
This is immediate for characteristic functions and follows for
$L^1(\mu)$ functions by the usual approximation process. We will
always assume that 
% $\mu(\widehat{\R})=0$, equivalently
$\mu(X')=\mu(\inte{X'})$ for any $j$-tile $X'$.  
The (measure theoretic)
\emph{entropy}\index{entropy}\label{not:entropy} of $\mu$ is then
given by 
\begin{equation}\label{eq:defentropy}
  h=h_{\mu}=\lim_{j\to\infty}-\frac{1}{j}\sum_{X'_j\in\X'_j}\mu(X'_j)\log\mu(X'_j).
\end{equation}
It is not very hard to show that the limit exists and is nonnegative. 
If $\mu$ is ergodic the Shannon--McMillan--Breiman
\index{Shannon-McMillan-Breiman} 
theorem (see, for example \cite{Petersen}) says that
\begin{equation}\label{eq:ShannMcM}
  -\frac{1}{j}\log\mu(X'_j)\to h, \text{ as } j\to\infty,
\end{equation}
 for $\mu$-almost every $\{x\}=\bigcap X'_j$, where $X'_j\in\X'_j$.
% and $X'_0\supset X'_1\supset\dots$. 

%% The \emph{Jacobien} of the measure $\mu$ is given by 
%% \begin{equation}\label{eq:JacDef}
%%   J_\mu(x):=\frac{d\mu\circ R}{d\mu}=\lim_j \frac{\mu (R(X'_j))}{\mu (X'_j)}.
%% \end{equation}
%% Here the first term denotes the Radon-Nikodym derivative, in the second term $x\in X'_j\in\X'_j$. The Jacobien describes how $R$ expands with respect to $\mu$. We have
%% \begin{equation}
%%   \label{eq:Jdmu}
%%   \int J_{\mu}\,d\mu=\int d\mu\circ R=\deg R.
%% \end{equation}
%%  The following is known as \emph{Rohlin's formula}:
%% \begin{equation}\label{eq:rohlin}
%%   h_\mu=\int\log J_\mu\,d\mu.
%% \end{equation}
%% By \emph{Jensen's inequality}
%% \begin{equation}\label{eq:Jensen}
%%   \int\log J_\mu\,d\mu\leq\log\int J_{\mu}\,d\mu=\log\deg R,
%% \end{equation}
%% where equality occurs if and only if $J_{\mu}=\const=\log R$. The measure satisfying this is called the measure of \emph{maximal entropy} (it is actually the image of Hausdorff measure on the snowsphere---in the Hausdorff dimension of it). 

The \emph{spherical derivative}\index{spherical derivative} of $R$ is given by
\begin{equation*}\label{not:sphericder}
  R^{\#}(z):=\abs{R'(z)}\frac{1+\abs{z}^2}{1+\abs{R(z)}^2},
\end{equation*}
here (and only here) the right-hand side is expressed in terms of the
Euclidean metric.   
It satisfies the usual rules (chain rule, derivative of inverse), as
well as 
\begin{align*}
  \left(\frac{1}{R}\right)^{\#}&=R^{\#} \text{ and}
  &R\left(\frac{1}{z}\right)^{\#}&=R^{\#}\left(\frac{1}{z}\right).
\end{align*}
The area element on the sphere satisfies
\begin{equation}
  \label{eq:DRw}
  d\lambda(R(w))=(R^{\#})^2d\lambda(w).   
\end{equation}

The \emph{Lyapunov exponent}\index{Lyapunov exponent} measures
metrical expansion,  in $(S^2,\abs{x-y})$ it is
\begin{align}\label{eq:defLyap}
  &\chi=\chi_\mu:=\int\log R^{\#}\,d\mu,
  \intertext{in $(S^2,\abs{x-y}_{\SC})$ it is}
  &\chi_{\SC}=\chi_{\mu,\SC}=\log N, \quad \text{by Theorem \ref{thm:dS2}
    (\ref{item:dS2Rsimil}).} 
%   \\
%   &=\int \log \abs{R'}\,d\mu+\int
%   \log(1+\abs{z}^2)\,d\mu-\int\log(1+\abs{R(z)}^2)\,d\mu \notag 
%   \\
%   &=\int \log \abs{R'}\,d\mu.
\end{align}
The Lyapunov exponent $\chi$
is nonnegative (see \cite{przy}). 
If $R$ is ergodic, we obtain by Birkhoff's ergodic theorem and the chain rule
\begin{equation}\label{eq:chiBirk}
  \frac{1}{k}\log(R^k)^{\#}(z)=\frac{1}{k}\sum_{j=0}^{k-1}\log
  R^{\#}(R^jz)\to \chi, \text{ as } k\to\infty 
\end{equation}
for $\mu$-almost every $z\in S^2$. Essentially by combining equations
(\ref{eq:ShannMcM}) and (\ref{eq:chiBirk}) one obtains Manning's
formula 
\begin{equation}\label{eq:manning}
  \dim \mu|(S^2,\abs{x-y}_{\SC}) =\frac{h}{\chi_{\SC}},
\end{equation}
when $\mu$ is ergodic and $\chi_\mu>0$. This was originally proved
(under stronger assumptions) in \cite{Manning},  the general case was
done in \cite{mane}. A complete proof can also be found in
\cite{przyurban}, Chapter 10.4.
% (the numbering may still change
%though). 
We will not directly use this formula.

\medskip
Since the critical values of $R$ are repelling fixed points, the Julia
set of $R$ is the whole sphere $\CDach$ (by the classification of
Fatou components, see \cite{Carleson}). The map $R$ then is
\emph{ergodic} (see Theorem 3.9 of \cite{McM}). 

\section{The Invariant Measure on the Sphere}\index{invariant measure}
\label{sec:invar-meas-sphere}

Lebesgue measure on the sphere is not invariant under the rational map
$R$. There is however an invariant measure $\mu$ that is absolutely
continuous with respect to ($2$-dimensional) Lebesgue measure
$\lambda$ on the sphere. 
% Since $\dim\mu=\dim\lambda=2$, we
% have  
% \begin{equation}
%   \label{eq:3hchi}
%   h_{\mu}=2\chi_{\mu}  
% \end{equation}
% by Manning's formula (\ref{eq:manning}). 
% By equations (\ref{eq:hmuS2}) (as well as Manning's formula) the
% dimension of $\mu$ in $(S^2,d')$ is
% \begin{equation}
%   \label{eq:dimmuS2}
%   \dim \mu|(S^2,d')= \frac{h_\mu}{\log N}=\frac{2 \chi_{\mu}}{\log N}.
% \end{equation}
% Recall that this is equal to the dimension of the elliptic harmonic
% measure of the corresponding snowsphere $\SC$. To prove Theorem
% \ref{thm:main} is it therefore enough to show that $h_{\mu}=2\log
% \chi_{\mu}< 2\log N$ by Theorem \ref{thm:dS2} (\ref{item:HdimS2}). 
% Equation (\ref{eq:fftilde}) implies that the
% dimension of the elliptic harmonic measure obtained from $f\colon
% \SC\to S^2$ is the same as the dimension of the measure $\mu:=\mu\circ
% \tilde{f}$, where $\tilde{f}\colon\widetilde{\SC}\to\Hob$. To not
% introduce more terminology $\mu$ will be called an elliptic harmonic
% measure as well. 
% Since $R$ commutes with the maps $\sigma_k\colon \widetilde{X}_k\to
% \widetilde{\SC}$ (encoding the self similarity, see section
% \ref{sec:genrat}), we have 
% \begin{equation*}
%   h_{\mu}=h_{\mu}.
% \end{equation*}
% The Lyapunov exponent of $\mu$ is known (by equation (\ref{eq:lyap2})). By Manning's formula it is therefore enough to calculate the entropy $h_{\mu}$, or equivalently $\chi_{\mu}$, to get the dimension of the elliptic harmonic measure. 
We now proceed to define this measure.
The technique to construct $\mu$ goes back to Sullivan and 
Patterson, the author learned it from \cite{grcsmir}. 
\begin{remark}
  In the following, $R^{-1}(z)$ will denote two different things:  the
  set of preimages, and a branch of the inverse function. We write
  $\{R^{-1}(z)\}$ to denote the former, $R^{-1}(z)$ for the latter. 
\end{remark}
For every $j\geq0$, let $d\mu_j:=\kappa_j\,d\lambda$, where
\begin{equation*}\label{not:kappaj}
  \kappa_j(z):=\sum_{w\in\{R^{-j}(z)\}}\left[(R^j)^{\#}(w)\right]^{-2},
\end{equation*}
and $\lambda$ is \emph{normalized} Lebesgue measure on the sphere $S^2$.
Note that by the chain rule (for $j\geq 1$)
\begin{equation}
  \label{eq:kjkj-1}
  \kappa_j(z)=\sum_{y\in \{R^{-1}(z)\}} R^{\#}(y)^{-2}\kappa_{j-1}(y). 
\end{equation}

\begin{lemma}
  \label{lem:intduk1}
  The measures $\mu_j$ are probability measures and satisfy
  \begin{align*}
    R_*\mu_{j-1}=\mu_{j},
  \end{align*}
  for every $j\geq 1$.
\end{lemma}
\begin{proof}

%   To see the first equality we compute
%   \begin{align*}
%     \int
%     d\mu_j&=\int\kappa_k(z)\,d\lambda(z)=\int\sum_{w\in\{R^{-j}(z)\}}((R^j)^{\#}(w))^{-2}\,d\lambda(z)     
%     \\
%     &=\sum_{X'\in\X'_j}\int_{X'}((R^j)^{\#}(w))^{-2}\,d\lambda(R^jw)
%     \\
%     &=\sum_{X'\in\X'_j}\int_{X'}((R^j)^{\#}(w))^{-2}((R^j)^{\#}(w))^{2}\,d\lambda(w)=1. 
%   \end{align*}
%   Where equation (\ref{eq:DRw}) was used.
%  Here property \ref{lem:sdl} of lemma \ref{lem:sphericalderprop} is used.
%   For the second equality let $A\subset S^2$ be a Borel set, then
%   \begin{align*}
%     \mu_{j+1}(A)&=\int_A\sum_{w\in \{R^{-j-1}z\}}((R^{j+1})^{\#}(w))^{-2}d\lambda(z)
%       \\
%       &=\int_A\sum_{z'\in \{R^{-1}z\}}(R^{\#}(z'))^{-2}\biggl(\sum_{w\in\{R^{-j}z'\}}((R^{j})^{\#}(w))^{-2}\biggr)d\lambda(z)
%       \\
%       &\overset{Rz'=z}{=}\sum_{B\in\{R^{-1}A\}}\int_B\sum_{w\in\{R^{-j}z'\}}((R^{j})^{\#}(w))^{-2}d\lambda(z')
%       \\
%       &=\mu_j(R^{-1}(A))=R_*\mu_j(A).
%   \end{align*}
  By equation (\ref{eq:kjkj-1}),
  \begin{align*}
    d\mu_j(z) & = \kappa_j(z)d\lambda(z)=\sum_{y\in \{R^{-1}(z)\}}
      R^{\#}(y)^{-2}\kappa_{j-1}(y) d\lambda(z),
      \intertext{since $z= R(y)$, we have  
      $d\lambda(z)= R^{\#}(y)^2d\lambda(y)$,} 
      \\
      & = \sum_{y\in \{R^{-1}(z)\}}\kappa_{j-1}(y) d\lambda(y).
  \end{align*}
  Thus, for any Borel set $A$
  \begin{equation*}
    \mu_j(A)=\int_A d\mu_j=\int_{\{R^{-1}A\}} \kappa_{j-1}d\lambda=
      \mu_{j-1}(R^{-1} A)= R_*\mu_{j-1}(A). 
  \end{equation*}
  A pushforward of a probability measure is a probability measure,
  thus all $\mu_j$ are probability measures.
\end{proof}
So by iteration $$\mu_j=R^j_*\lambda,$$ which we could have used as the definition.
Since the set of probability measures is compact in the weak-$*$ topology, we can define $\mu$ as a weak-$*$ subsequential limit of
\begin{equation}
  \label{eq:mujav}
  \overline{\mu}_j=\frac{1}{j}(\mu_1+\dots+\mu_j).
\end{equation}
Since 
\begin{align*}
  \left|\frac{1}{j}(\mu_1+\dots+\mu_j)(R^{-1}A)-\frac{1}{j}(\mu_1+\dots+\mu_j)(A)\right|
    \\
    \leq
    \frac{1}{j}\mu_{j+1}(A)+\frac{1}{j}\mu_1(A)
    \leq\frac{2}{j},
\end{align*}
for any Borel set $A\subset S^2$, the measure $\mu$ is invariant. It
is a probability measure as the weak-$*$ limit of probability measures
on a compact space.
We need 
to show that $\mu$ is absolutely continuous with respect to
Lebesgue measure.
In the next lemma, it is shown that there is $G\in L^1(\lambda)$
such that $\kappa_j\leq G$.  This shows absolute continuity of $\mu$.
%as well as that $\mu$ is a probability measure (by the dominated
%convergence theorem).  

% In
% the next lemma we will see that the $\kappa_j$ converge locally
% uniformly (away from the postcritical set).  
%It is however not yet clear that $\mu$ is absolutely continuous with respect to Lebesgue measure $\lambda$. 
% We will also show that there is a function $G(z)\colon \CDach\to
% (0,\infty]$, with $\kappa_j\leq G$ and $G\in L^1(\lambda)$. From this
% follows that $\mu$ is absolutely continuous with respect to Lebesgue
% measure $\lambda$ (on $S^2$), by the dominated convergence theorem. 

\medskip
To estimate the $\kappa_j$, we need to take preimages, more precisely
estimate derivatives along inverse orbits. Since the rational map $R$
is postcritically finite, this is not too hard.  
A branch of the inverse $R^{-1}$ is defined in a (simply connected)
neighborhood $U$ of a point $z$ if and only if $U$ does not contain a
critical value. 
Every preimage of such a $U$ will again not contain any postcritical point
by the postcritical finiteness. 
%For a point $z$ is that is not a postcritical point, all of the preimages will also not be postcritical. In fact for any non-postcritical $z$ (i.~e.\ $z\neq1,-1,\infty$), there in a neighborhood $U_{\epsilon}(z)$ of $z$ containing no postcritical point. 
Thus, every branch of $R^{-j}$ is then (univalently) defined on
$U$. We can then use Koebe's theorem to control the derivative.   

% In the following ``$\diam$'', ``$\dist$'', and ``$\area$''  refer to
% the spherical metric on $\CDach$, $\abs{x-y}$ is the Euclidean metric
% on $\C$. The $\epsilon$-neighborhoods $U_{\epsilon}$ are defined with
% respect to the spherical metric. 
 
By $\deg(c)=\deg_R(c)$, we denote the \emph{degree} of a critical point
$c$, i.e.,  
\begin{equation*}\label{eq:Rexp}
  R(w)=p+a(w-c)^{\deg(c)}+\dots,
\end{equation*}
in a neighborhood of $c$. It will simplify our discussion to allow
$\deg(c)=1$. 
For a postcritical point $p\in\post(R)$, 
let 
\begin{equation*}\label{not:maxdeg}
  \max\deg(p):=\max\{\deg_{R^j}(c):c\in\{R^{-j}(p)\}, j\geq 1\}.
\end{equation*}
This is finite, since $R$ has no critical periodic orbits (see Lemma
\ref{lem:propR}).

\begin{lemma}
  \label{lem:densitybehavior} 
  The densities $\kappa_j\colon\CDach\to
  \R$ satisfy the following. 
  \begin{enumerate}
  \item 
    \label{lem:equi} 
    They are \emph{equicontinuous} on
    $\CDach\setminus\post(R)$. 
  \item 
    \label{lem:sasymp1}
    Let $z\in \CDach\setminus \post(R)$, $\dist(z, \post(R))\geq \epsilon$.
    %$U=\bigcup_{p\in\post(R)}U_{\epsilon}(p)$ be an
    %$\epsilon$-neighborhood of the postcritical set. 
    Then  
    \begin{equation*}
      \kappa_j(z)\asymp 1,      
    \end{equation*}
    %for $z\in\CDach\setminus U$, 
    where $C(\asymp)=C(\epsilon)$ is
    independent of $j$.  
  \item \label{lem:sasymppole}
    In a neighborhood of a postcritical point $p$, 
    \begin{equation*}
      \kappa_j(z)\asymp\abs{z-p}^{-2\left(1-\frac{1}{\max\deg(p)}\right)},
    \end{equation*}
    for $j\geq k_0$, with $C(\asymp)$ independent of $j$.
  \end{enumerate}
\end{lemma}
\begin{proof}
% (1)
(\ref{lem:equi})  
Consider an arbitrary $z\in\CDach\setminus\post(R)$. Let
$\delta:=\dist(z,\post(R))$. Every branch of $R^{-j}$ is then defined on
the neighborhood $U_\delta(z)$. 
Let $\epsilon< \delta$, then for $z'\in
U_{\epsilon}(z)\subset U_{\delta}(z)$ we have by Koebe distortion,

\begin{equation*}
  ((R^j)^{\#}(w))^{-1}=(R^{-j})^{\#}(z)\asymp(R^{-j})^{\#}(z')=
  ((R^j)^{\#}(w'))^{-1},   
\end{equation*}
where $w=R^{-j}(z)$, $w'=R^{-j}(z')$.
Here $C(\asymp)=C(\epsilon/\delta)\to 1$ as  $\epsilon/\delta\to 0$
(independent of $j$). 
This yields
\begin{equation*}
  \kappa_j(z)=\sum_{w\in\{R^{-j}(z)\}}((R^j)^{\#}(w))^{-2}\asymp
  \sum_{w'\in\{R^{-j}(z')\}}((R^j)^{\#}(w'))^{-2}=\kappa_j(z'), 
\end{equation*}
(with $C(\asymp)=C(\epsilon/\delta)\to 1$ as $\epsilon/\delta\to 0$ as before).
% (independent of $j$). 
The statement follows from (\ref{lem:sasymp1}) which we prove next.

\medskip
% (2)
(\ref{lem:sasymp1})
Let $z\in \CDach \setminus \post(R)$. Let $Y'_1,Y'_2$ be the two
$0$-tiles, they are $\Hobd,\Hubd$ if $R$ is as in Section
\ref{sec:rati-map-repr}. Without loss of generality, $z\in Y'_1$. 
 
Assume first that $z\in Z'_1\in\X'_1$, where the   $1$-tile $Z'_1$
does not contain a postcritical point. Let $U$ be a neighborhood of
$Z'_1$ containing no postcritical point. Fix a branch of $R^{-j}$ on
$U$. Let $w:=R^{-j}(z)$ and $R^{-j}(Z'_1)=:Z'_{j+1}\in \X'_{j+1}$. 
There is (exactly) one $j$-tile $X'_j\ni w$ such that
$R^j(X'_j)=Y'_1$; conversely for each $j$-tile $X'_j$ with
$R^j(X'_j)=Y'_1$ there is (exactly) one such $w\in X'_j$. 
By
Koebe distortion 
%and property \ref{lem:sdinv} of lemma \ref{lem:sphericalderprop}
\begin{align}
  \notag
  \diam Z'_{j+1}&\asymp (R^{-j})^{\#}(z)\diam Z'_1\asymp ((R^j)^{\#}(w))^{-1}
  \intertext{ and by Lemma \ref{lem:combmetric} }
  \label{eq:prooflemmakappa1}
  \area X'_j & \asymp \area Z'_{j+1} \asymp((R^j)^{\#}(w))^{-2}.
\end{align}
%for a $j$-tile $X'_j\cap Z'_{j+1}\neq \emptyset$ with
%$R^j(X'_j)=Y'_1$. 
Here, $C(\asymp)$ is independent of $j$ or the branch of $R^{-j}$. 

\smallskip
By the above, we can estimate the area of all $j$-tiles $X'_j$ that are
preimages of $Y'_1$. To estimate the area of $j$-tiles
$\widetilde{X}'_j$ that are
preimages of $Y'_2$, we \defn{match} each such $j$-tile
$\widetilde{X}_j$ to one $X'_j$ (this is in fact a \defn{perfect
  matching}).   
Fix a $0$-edge $E$ (which is $[-\infty, -1]$, $[-1,1]$, or
$[1,\infty]$ for $R$ as in Section \ref{sec:rati-map-repr}). Consider
now a $j$-tile $\widetilde{X}'_j$ that is mapped to $Y'_2$,
$R^j(\widetilde{X}'_j)= Y'_2$. Then there is exactly one 
$j$-tile $X'_j$ (which is mapped to $Y'_1$ by $R^j$) that shares a
preimage of $E$ with $\widetilde{X}'_j$. By Lemma \ref{lem:combmetric},
 $\area \widetilde{X}'_j \asymp \area X'_j$. Thus,  
% Assume wlog that $Z'_1\subset \Hobd$. 
% Consider a $j$-tile $X'_j$ such that $R^j(X'_j)=\Hobd$. 
% Then every $(j+1)$-tile $Z'_{j+1}$ as above is contained in one such
% $j$-tile $X'_j$. Thus $\diam X'_j\asymp \diam Z'_{j+1}$. 

% Let $E_j$ be a $j$-edge such that $R^j(E_j)= [-1,1]$. 
% Consider now a $j$-tile $Y'_j$ such that $R^j(Y''_j)=\Hubd$. Then there
% is exactly one $j$-tile $X'_j$ (with $R^j(X'_j)=\Hobd$) sharing $E_j$
% with $Y''_j$. Recall that $\diam Y''_j\asymp \diam X'_j$ (ref needed). 
% Thus
\begin{align*}
  1&=\sum_{X'_j\in \X'_j}\area X'_j
  \asymp\sum_{\substack{X'_j\in \X'_j\\ R^j(X'_j)=Y'_1}}\area X'_j, \quad
  \text{ thus by (\ref{eq:prooflemmakappa1})}
  % \asymp\sum_{Z'_{j+1}\in\{R^{-j}(Z'_1)\}}\area Z'_{j+1}
  \\
  &\asymp \sum_{w\in \{R^{-j}(z)\}}((R^j)^{\#}(w))^{-2}=\kappa_j(z).
\end{align*}
To complete the statement, we consider $z\in Z'_k\in\X'_k$, where the
$k$-tile $Z'_k$ does not contain a postcritical point, and repeat the
argument. Every point $z\in \CDach\setminus \post(R)$ with
$\dist(z,\post(R))\geq \epsilon$ will be contained in such a $Z'_k$
for $k=k(\epsilon)$ sufficiently large by Lemma
\ref{lem:jtilesshrink}.

\bigskip
(\ref{lem:sasymppole})
% 3
Each postcritical point has a finite orbit, since $R$ is postcritically
finite. Let $\postpre(R)$ be the set of preperiodic postcritical
points of $R$, and $\postper(R)$ be the set of periodic postcritical
points of $R$,
\begin{equation*}
  \post(R)=\postpre(R)\cup \postper(R).
\end{equation*}
Each point $p\in \postper(R)$ is a (repelling) fixed point for a
suitable iterate $R^{n_0}$.  
The proof for (\ref{lem:sasymppole}) is broken up into three parts. In
(\ref{lem:sasymppole}b), we assume that
each $p\in \postper(R)$ is a fixed point. 

\medskip
% 3a
(\ref{lem:sasymppole}a) 
Consider first a (preperiodic) point
$p\in\postpre(R)$. Then there is a $k_0\geq 1$ such that no $c\in
\{R^{-k_0}(p)\}$ is a postcritical point. Let $z\in U_\epsilon(p)$,
where $\epsilon$ is smaller than the distance of $p$ to another
postcritical point. 
Then for any $y\in \{R^{-k_0}(z)\}$ 
\begin{align*}
  & \abs{z-p}=\abs{R^{k_0}(y)-p}\asymp \abs{y-c}^d, \quad \text{and}
  \\
  & (R^{k_0})^{\#}(y)\asymp \abs{y-c}^{d-1}\asymp \abs{z-p}^{\frac{d-1}{d}},
\end{align*}
where $c\in \{R^{-k_0}(p)\}$, $d=d(y)=\deg_{R^{k_0}}(c)$, and $C(\asymp)$
does only depend on $\epsilon$. Thus, (see (\ref{eq:kjkj-1}))

\begin{align*}
  \kappa_{j+k_0}(z)& = \sum_{y\in \{R^{-k_0}(z)\}}
    (R^{k_0})^{\#}(y)^{-2}\underbrace{\kappa_j(y)}_{\asymp 1 \text{ by
      (\ref{lem:sasymp1})}}
  \\
  & \asymp \sum_{y\in
    \{R^{k_0}(z)\}}\abs{z-p}^{-2\left(1-\frac{1}{d(y)}\right)}
  \asymp \abs{z-p}^{-2\left(1-\frac{1}{\max\deg(p)}\right)}.
\end{align*}
For $j< k_0$,
the same argument as above yields  $\kappa_j(z) \lesssim
\abs{z-p}^{-2\left(1-\frac{1}{\max\deg(p)}\right)}$.

\medskip
% 3b
(\ref{lem:sasymppole}b)
Consider now a  point $p\in \postper(R)$. 
We assume here that $p$ is a (repelling) fixed point of $R$.
The
\emph{multiplier} (of $R$ at $p$) is defined as  
\begin{equation*}
  \Lambda=\Lambda(p):=\abs{R'(p)}=R^{\#}(p)> 1.
\end{equation*}
By Koenigs' linearization theorem (see \cite{Carleson}, II.2 and II.3),
there is a neighborhood $U=U(p)$ of $p$ in which $R$ is conformally
conjugate to $z\mapsto R'(p)z$. Let $R^{-k}_p$ be the branch of
$R^{-k}$ that keeps $p$ fixed (defined on $U$).  

Let $z\in U$, and set $z_k:=R^{-k}_p(z)$. Then (by Koebe) 
\begin{equation*}
  \abs{z_k-p}\asymp \Lambda^{-k}\abs{z-p} \quad \text{and }
  (R^k)^{\#}(z_k)\asymp \Lambda^k,
\end{equation*}
where $C(\asymp)$ is independent of $k$. Consider now a point $w\in
\{R^{-1}(z_k)\}$, $w\neq z_{k+1}$. Then $w$ is in the neighborhood of
a point $p'\in \{R^{-1}(p)\}$, where $p'\notin \postper(R)$, and with
$d=d(w)=\deg_R(p')$ 
\begin{align}  
  \notag
  & \abs{w-p'}^d \asymp \abs{z_k-p}\asymp \Lambda^{-k}\abs{z-p},
  \\
  \notag
  & R^{\#}(w) \asymp \abs{w-p'}^{d-1}\asymp
  \Lambda^{-k\left(\frac{d-1}{d}\right)}\abs{z-p}^{1-\frac{1}{d}},
  \\
  \notag
  & (R^{k+1})^{\#}(w)  =R^{\#}(w) (R^k)^{\#}(z_k)  
  \asymp \Lambda^{\frac{k}{d}}\abs{z-p}^{1-\frac{1}{d}}.
  \intertext{Thus by (3a) for $j\geq k_0$} 
  \notag
  & \kappa_j(w)  \asymp
  \abs{w-p'}^{-2\left(1-\frac{1}{\max\deg(p')}\right)}
  \\
  \notag
  & \phantom{XXX} \asymp
  \Lambda^{\frac{2k}{d}\left(1-\frac{1}{\max\deg(p')}\right)}
  \abs{z-p}^{\frac{-2}{d}\left(1-\frac{1}{\max\deg(p')}\right)}
  \\
  \label{eq:kappajestimate}
  & \left[(R^{k+1})^{\#}(w)\right]^{-2}\kappa_j(w) 
  \asymp \Lambda^{\frac{-2k}{d\max\deg(p')}}
    \abs{z-p}^{-2\left(1-\frac{1}{d\max\deg(p')}\right)}. 
\end{align}
Here, $C(\asymp)$ does only depend on $U$. 
If $j<k_0$, we obtain in the previous two expressions ``$\lesssim$''
instead of ``$\asymp$''.
Thus, (using (\ref{eq:kjkj-1}))
\begin{align*}
  \kappa_j(z)= & \sum_{\substack{w\in\{R^{-1}(z)\}\\ w\neq
      z_1}}R^{\#}(w)^{-2}{\kappa_{j-1}(w)} 
  \\
  &+\sum_{\substack{w\in\{R^{-1}(z_1)\}\\ w\neq
      z_2}}\left[(R^{2})^{\#}(w)\right]^{-2}{\kappa_{j-2}(w)} 
  \\
  &\vdots
  \\
  &+\sum_{\substack{w\in\{R^{-1}(z_{j-1})\}\\
      w\neq z_{j}}}\left[(R^{j})^{\#}(w)\right]^{-2}\underbrace{\kappa_{j-j}(w)}_{=
    1} 
  \\
  &+\left[(R^j)^{\#}\big(z_{j}\big)\right]^{-2}.
  \\
  \intertext{ Hence by (\ref{eq:kappajestimate}) for $j> k_0$}
  \asymp & \abs{z-p}^{-2\left(1-\frac{1}{\max\deg(p)}\right)} 
      \underbrace{\left(1+ \Lambda^{\frac{-2}{\max\deg(p)}} + \dots +
          \Lambda^{\frac{-2(j-1)}{\max\deg(p)}} \right)}_{\asymp 1}
        + \Lambda^{-2j}.
\end{align*}
This proves property (\ref{lem:sasymppole}) in this case.

\medskip
% 3c
(\ref{lem:sasymppole}c) 
Consider now a $p\in \postper(R)$, where 
we now drop the assumption that $p$ is a fixed point $R$. 
Let $n_0$ be the period of $p$. Then we know 
the behavior of $\kappa_{jn_0}$
near $p$ by (\ref{lem:sasymppole}b). 
It is therefore enough
to show that if $\kappa_j$ has the desired behavior near $p$, then
$\kappa_{j+1}$ has this behavior as well. 

Let $z\in U(p)$.  
Then $\{R^{-1}(z)\}$ contains one point $z'$ in a neighborhood of $p'\in
\postper(R)$. We have $R^{\#}(z')\asymp 1$, and
$\max\deg(p)=\max\deg(p')$. 
Thus, 
\begin{align*}
  & \abs{z'-p'}\asymp \abs{z-p}.
  \\
  \intertext{We assume now that for $j\geq k_0$}
  & \kappa_j(z')\asymp \abs{z'-p'}^{-2\left(1-\frac{1}{\max \deg
        p'}\right)},
  \text{ this yields}
  \\
  & \phantom{\kappa_j(z')} \asymp \abs{z-p}^{-2\left(1-\frac{1}{\max \deg p}\right)}.
\end{align*}

Consider now $w\in
\{R^{-1}(z)\}, w\neq z'$. Note that $w$ lies in a neighborhood of $q\in
\{R^{-1}(p)\}$, where $q\notin \postper(R)$. Then we obtain as before with
$d=\deg(q)$ and a constant $C(\asymp)$ (only depending on $U(p)$)
\begin{align*}
  & \abs{w-q}^d \asymp\abs{z-p}, 
  \\
  & R^{\#}(w) \asymp \abs{w-q}^{d-1}\asymp \abs{z-p}^{1-\frac{1}{d}},
  \\
  & \kappa_j(w) \asymp
  \abs{w-q}^{-2\left(1-\frac{1}{\max\deg(q)}\right)},
  \text{ from (\ref{lem:sasymppole}a) for $j\geq k_0$,}
  \\
  & \phantom{\kappa_j(w)} \asymp
  \abs{z-p}^{-\frac{2}{d}\left(1-\frac{1}{\max\deg(q)}\right)},  
  \\
  &R^{\#}(w)^{-2}\kappa_j(w) \asymp
  \abs{z-p}^{-2\left(1-\frac{1}{d\max\deg(q)}\right)}.  
\end{align*}
Thus, (using (\ref{eq:kjkj-1})) for $j\geq k_0$
\begin{align*}
  \kappa_{j+1}(z) & =R^{\#}(z')^{-2}\kappa_j(z') 
  + \sum_{\substack{w\in \{R^{-1}(z)\} \\ w\neq z'}}
  R^{\#}(w)^{-2}\kappa_j(w)
  \\
  & \asymp \abs{z-p}^{-2\left(1-\frac{1}{\max\deg{p}}\right)}.
\end{align*}

\end{proof}

\begin{lemma}
  \label{lem:kjtok}
  The averages $\overline{\mu}_j:=\frac{1}{j}\left(\mu_1+\dots +\mu_j\right)$ 
%of the measures $\mu_j$ (see (\ref{eq:mujav})) 
  converge 
  weak-$*$ to an $R$-invariant ergodic probability measure $\mu$.
  Furthermore,
  \begin{enumerate}
  \item $d\mu=\kappa d\lambda$ is
    absolutely continuous with respect to Lebesgue measure.
  \item 
    \label{lemkasympcont}
    The density $\kappa$ is continuous on $\CDach\setminus \post(R)$. 
  \item \label{lem:kasymp1}
    Let
    $U=\bigcup_{p\in\post(R)}U_{\epsilon}(p)$ be an
    $\epsilon$-neighborhood of the postcritical set. Then  
    \begin{equation*}
      \kappa(z)\asymp 1,      
    \end{equation*}
    for $z\in\CDach\setminus U$, where $C(\asymp)=C(\epsilon)$. 
  \item \label{lem:kasymppole}
    In a neighborhood of a postcritical point $p$, 
    \begin{equation*}
      \kappa(z)\asymp\abs{z-p}^{-2\left(1-\frac{1}{\max\deg(p)}\right)}.
    \end{equation*}
  \end{enumerate}
\end{lemma}

\begin{proof}
  Every weak-$*$ subsequential limit $\mu$ of $(\overline{\mu}_j)$
  satisfies (\ref{lem:kasymp1}) and   
  (\ref{lem:kasymppole}) by the last lemma. Note that all $\kappa_j$
  are dominated by a 
  $G\in L^1(\lambda)$. Thus, $\mu$ is absolutely continuous with respect to
  Lebesgue measure. 
  It is well known (see \cite{McM}, Theorem 3.9) that $\mu$ is
  ergodic, since 
  the Julia set of $R$ is $\CDach$. Recall that two ergodic invariant
  probability 
  measures are either singular or identical. Thus, $(\overline{\mu}_j)$
  converges. The measure $\mu$ is a probability measure as a
  weak-$*$ limit of probability measures on a compact space. 
  %by the
  %dominated convergence theorem.

  The averages $\frac{1}{j}(\kappa_1+\dots +\kappa_j)$ are equicontinuous on
  $\CDach\setminus\post(R)$ by the last lemma. Taking another
  subsequence from the above, we obtain by Arzel\`{a}-Ascoli that
  $\kappa$ is continuous. Uniqueness of ergodic, nonsingular measures
  yields (\ref{lemkasympcont}).  
\end{proof} 

\section{Proof of the Theorems}
\label{sec:estim-hold-expon}
In this section, we prove Theorem \ref{thm:Ralpha} and  Theorem
\ref{thm:Rmu} (and thus, Theorem \ref{thm:almostbiHoelder} and Theorem
\ref{thm:main}) for 
\begin{equation*}
  \alpha= \frac{\log N}{\chi}.
\end{equation*}
Recall that $\SC$ was constructed from an $N$-generator. 
%In the next section $\alpha > \frac{2}{\dim_H(\SC)}$ is shown. 
Let us first express
the entropy and the Lyapunov exponent in terms of the invariant
measure $\mu$ constructed in the previous section.

\smallskip
The \emph{Jacobian}\index{Jacobian} describes the expansion of $R$
with respect to $\mu$, 
\begin{equation}
  \label{eq:2defJac}
  J_{\mu}(x):=\frac{d\mu\circ R}{d\mu}=\lim_j \frac{\mu
    (R(X'_j))}{\mu(X'_j)}=R^{\#}(z)^2\frac{\kappa(Rz)}{\kappa(z)}.
\end{equation}
Here, the first term denotes the Radon--Nikodym derivative, in the
second term $x\in X'_j\in\X'_j$. 
Note that
\begin{equation}
  \label{eq:Jdmu}
  \int J_{\mu}\,d\mu=\int d\mu\circ R=\sum_{X'\in
    \X'_1}\int_{X'}d\mu\circ R      \overset{Rz=w}{=}\deg R. 
\end{equation}
The entropy may be expressed via the Jacobian
\begin{equation}
  \label{eq:Rohlinsform}
  h=h_\mu = \int \log J_\mu \,d\mu. 
\end{equation}
This is known as \emph{Rohlin's formula} (see for example \cite{przyurban}, Chapter 1).  
Thus, (using (\ref{eq:invmeas2}))  
\begin{align}
  \label{eq:hchi}
  h&=\int\log(R^{\#})^2\,d\mu+\int\log\kappa(Rz)\,d\mu-\int
  \log\kappa(z)\,d\mu = 2\chi.
\end{align}
This is a trivial instance of Manning's formula (\ref{eq:manning}), it
just says that 
Lebesgue measure (on $S^2$) has dimension $2$. 

By (\ref{eq:hchi}), (\ref{eq:Rohlinsform}), (\ref{eq:Jdmu}), and
\emph{Jensen's inequality}\index{Jensen's inequality} 
\begin{equation}\label{eq:Jensen}
  2\chi=h=
  \int\log J_{\mu}\,d\mu\leq\log\int J_{\mu}\,d\mu=\log\deg R,
\end{equation}
where equality occurs if and only if $J_{\mu}=\const=\deg R$ ($\mu$
almost everywhere). The
measure satisfying this is called the measure of \emph{maximal
  entropy}, it is actually the Hausdorff measure on the
snowsphere---in the Hausdorff dimension of it. 
It was shown in \cite{ZdunikParabolicmaxmeas} that Lebesgue measure is
not absolutely continuous to the measure of maximal entropy unless $R$
is a \defn{Latt\`{e}s map} (has \defn{parabolic orbifold}). Thus,
$2\chi< \log \deg R$, unless $R$ is Latt\`{e}s in which case $2\chi=
\log \deg R$. However, we can derive this relatively easily 
given the setup of the last section. This is done for the convenience
of the reader in the next section. 

Recall that $\dim_H(\SC)=\frac{\log \deg R}{\log N}$ (Theorem \ref{thm:dS2}
(\ref{item:HdimS2})). Hence, 
\begin{align*}
  &\alpha= \frac{\log N}{\chi} = \frac{2}{\dim_H(\SC)}, &&\text{ if $R$
    is Latt\`{e}s}
  \\
  &\alpha > \frac{2}{\dim_H(\SC)}, &&\text{ otherwise}
\end{align*}
as desired.

% \medskip
% Let us
% recall how distances may be measured in combinatorial terms. Let
% $x,y\in S^2$ be arbitrary. Then 
% \begin{equation}
%   \label{eq:defjxy}
%   j(x,y):= \min \{j \mid \text{ there exist disjoint }
%   j\text{-tiles } X_j\ni x, Y_j\ni y\}.
% \end{equation}
% This measures to what level of tiles one needs to descend to be able to
% distinguish $x$ and $y$. The following are from
% \cite{snowballquasiball} (Lemma 3.10, Lemma 2.4, and Equation (6.10)). 
% \begin{align}
%   \label{eq:dS2xy}
%   &\abs{x-y}\asymp \diam X'_j\\
%   \label{eq:dSxy}
%   &\abs{x-y}_{\SC}\asymp N^{-j},
% \end{align}
% here $j=j(x,y)$, $C(\asymp)=C(N)$, and $X'_j$ is a $j$-tile containing
% $x$. Furthermore (for $j$-tiles $X'_j,Y'_j$)
% \begin{align}
%   \label{eq:diamXY}
%   &\diam X'_j\asymp \diam Y'_j, \quad \text{ if } X'_j\cap Y'_j\neq
%   \emptyset, 
%   \\
%   \label{eq:diamXX}
%   &\diam X'_{j+1} \asymp \diam X'_j, \quad \text{ for any $(j+1)$-tile
%   } X'_{j+1}\subset X'_j,
%   \\
%   \notag
%   %\label{eq:diamXareaX}
%   &\area(X'_j)\asymp (\diam X'_j)^2,
% \end{align}
% with a constant $C(\asymp)$. The same estimates hold (by the same
% arguments) for the $j$-tiles generated by any rational map (which is
% postcritically finite with no periodic critical point).

Let $\dist(x, \post(R))> \epsilon>0$ and $x\in X'_j\in \X'_j$. By
Lemma \ref{lem:combmetric} and 
Lemma \ref{lem:kjtok} (\ref{lem:kasymp1}), we obtain for sufficiently
large $j$
\begin{equation}
  \mu(X'_j)\asymp (\diam X'_j)^2,
\end{equation}
where $C(\asymp)=C(\epsilon)$. Thus using Shannon--McMillan--Breiman
(\ref{eq:ShannMcM}) 
\begin{equation}
  \label{eq:hdiamXj}
  h=\lim_j -\frac{1}{j}\log \mu(X'_j)=2 \lim_j -\frac{1}{j} \log \diam X'_j,
\end{equation}
for Lebesgue almost every $\{x\}=\bigcap X'_j$. 
We are now ready to prove Theorem \ref{thm:Ralpha}.

\begin{proof}[Proof of Theorem \ref{thm:Ralpha}]
 In the following, we write $a= b\pm C$ if $b-C\leq a\leq b+C$. 
 
 \smallskip
 Consider a $\{x\}=\bigcap X'_j$ satisfying (\ref{eq:hdiamXj}). Let
 $j=j(x,y)$ (see (\ref{eq:defjxy})), then $y\to x \Leftrightarrow j\to
 \infty$. Thus, 
 \begin{align*}
   \lim_{y\to x}&\frac{\log\abs{x-y}_{\SC}}{\log\abs{x-y}}
   = \lim_j \frac{-j\left(\log N \pm {C}/{j}\right)}{\diam X'_j \pm C}
   \quad \text{ by  (\ref{eq:dS2xy}) and (\ref{eq:dSxy}).}
   \\
   \intertext{Thus using (\ref{eq:hdiamXj}) and (\ref{eq:hchi})}
   & = \lim_j \frac{\log N \pm C/j}{-\frac{1}{j}\diam X'_j \pm
     C/j}=\frac{\log N}{h/2}=\frac{\log N}{\chi}=\alpha.
 \end{align*}
 
\end{proof}

\begin{proof}[Proof of Theorem \ref{thm:Rmu}]
  This
  is an easy consequence of Theorem
  \ref{thm:Ralpha}. Indeed consider a homeomorphism $f\colon \X\to
  \mathbf{Y}$ between metric spaces, such that
  \begin{align*}
    \frac{\log\abs{f(x)-f(y)}}{\log \abs{x-y}}\leq \alpha+\epsilon,
    \\
    \frac{\log\abs{f(x)-f(y)}}{\log \abs{x-y}}\geq \alpha -\epsilon,
  \end{align*}
  for all $0<\abs{x-y}<1/n$. Then we obtain for the Hausdorff
  dimension of the spaces 
  $\dim \mathbf{Y} \geq \dim\X/(\alpha+\epsilon)$ and $\dim
  \mathbf{Y} \leq 
  \dim\X/(\alpha-\epsilon)$. This follows directly from the definition
  of Hausdorff 
  dimension. Thus, Theorem \ref{thm:Rmu} is proved, since for every
  $\epsilon>0$ we can exhaust
  a set of full measure where the above holds (by Theorem
  \ref{thm:Ralpha}). 
\end{proof}

\begin{proof}[Proof of Theorem \ref{thm:lipschitz_lattes}]
  Let $(S^2, \abs{x-y}_{\SC})$ be snowflake equivalent to $S^2$. In
  particular they are quasisymmetric. Thus, by
  \cite{BonkMeyer_MarkovThurston}, the expanding Thurston map $R$ is
  topologically conjugate to a rational map. So we can assume that $R$
  is in fact a rational map (which is postcritically finite and has no
  periodic critical points). 

  Let $\varphi \colon (S^2,\abs{x-y}_{\SC})\to S^2$ be the snowflake
  equivalence. This means that there is a $\beta>0$, such that
  \begin{equation*}
    \abs{\varphi(x)-\varphi(y)}^\beta\asymp \abs{x-y}_{\SC},
  \end{equation*}
  for all $x,y\in S^2$. 
  Clearly $\varphi$ changes the (Hausdorff) dimension by the factor
  $\beta$. Thus we 
  have $\beta=2/d$, where $d=\dim_H (S^2,\abs{x-y}_\SC)$. The map
  $\varphi$ maps $d$-dimensional Hausdorff measure to
  ($2$-dimensional) Lebesgue measure (up to multiplicative constants).  

  Assume now that $R$ is not a Latt\`{e}s map. Then by Theorem
  \ref{thm:Rmu} there is a set $A\subset S^2$ that has full ($2$-dimensional)
  Lebesgue measure $\lambda$, yet is a (Hausdorff) $d$-dimensional
  zero set in $(S^2,\abs{x-y}_\SC)$. 

  The map $\id\colon S^2 \to (S^2,\abs{x-y}_{\SC})$ is a quasisymmetry
  (see \cite{BonkMeyer_MarkovThurston} and \cite{snowballquasiball}).  
  We now get a contradiction since
  the composition
  \begin{equation*}
    S^2\xrightarrow{\id} (S^2,\abs{x-y}_\SC) \xrightarrow{\varphi} S^2,
  \end{equation*}
  is a quasiconformal map that maps $A$ to a zero set (with respect to
  Lebesgue measure $\lambda$); hence is not absolutely continuous.

  \medskip
  Assume now that $R$ is topologically conjugate to a Latt\`{e}s
  map. So assume that $R$ is a Latt\`{e}s map. Then $R$ is obtained as
  a quotient of a linear map by a \emph{wallpaper group}. 
  We can choose $\abs{x-y}_{\SC}$ to be the 
  projection of the Euclidean metric. Then $(S^2,\abs{x-y}_{\SC})$ is
  easily seen to be bi-Lipschitz to $S^2$. 
  
  To be more explicit, we remind the reader of
  Theorem 3.1
  in \cite{MR2348953}.  The map $R$ may be described as follows. 
  \begin{itemize}
  \item There is a flat torus $\T^2\cong \C/\Lambda$, where
    $\Lambda\subset \C$ is a lattice of rank $2$.
  \item Furthermore, there is an affine map $L= az+b$, $a,b\in \C$,
    $a\ne 0$; satisfying 
    $L\Lambda\subset \Lambda$ (it holds $\abs{a}^2=\deg R$).
  \item Then $R$ is conformally conjugate to 
    \begin{equation*}
      L/G_n\colon \T^2/G_n\to \T^2/G_n,
    \end{equation*}
    for some $n$. 
    Here, $G_n$ is the group of $n$-th roots of unity, acting on $\T^2$
    by rotation around a base point. 
  \end{itemize}

  The space $\T^2/G_n$ is topologically a sphere. In fact, in the flat
  metric inherited from the torus $\T^2$ the sphere $\T^2/G_n$ is
  isometric to 
  (see Section 4 in \cite{MR2348953}) the path metric on
  \begin{itemize}
  \item a tetrahedron if the signature of $R$ is $(2,2,2,2)$.
  \item two triangles glued together along their boundaries, in the
    other cases. 
  \end{itemize}
  Clearly, these spheres are bi-Lipschitz to the standard sphere. Note
  that $\id\colon (S^2,\abs{x-y}_{\SC})\to S^2$ is \emph{not}
  bi-Lipschitz. 
\end{proof}

\section{The Jacobian of the Invariant Measure}
\label{sec:jacob-invar-meas}
% The \emph{Jacobian}\index{Jacobian} describes the expansion of $R$
% with respect to $\mu'$, 
% \begin{equation}
%   \label{eq:2defJac}
%   J_{\mu'}(x):=\frac{d\mu'\circ R}{d\mu'}=\lim_j \frac{\mu'
%     (R(X'_j))}{\mu' (X'_j)}. 
% \end{equation}
% Here the first term denotes the Radon-Nikodym derivative, in the
% second term $x\in X'_j\in\X'_j$. 
% We have
% \begin{equation}
%   \label{eq:Jdmu}
%   \int J_{\mu'}\,d\mu'=\int d\mu'\circ R=\sum_{X'\in
%     \X'_1}\int_{X'}d\mu'\circ R      \overset{Rz=w}{=}\deg R. 
% \end{equation}
%   By definition of the Jacobian and property \ref{lem:sdl} of lemma
%   \ref{lem:sphericalderprop} 
%   \begin{equation}\label{eq:Jeq}
%     J_{\mu'}(z)=(R^{\#}(z))^2\frac{\kappa(Rz)}{\kappa(z)}.
%   \end{equation}
% Therefore by equations (\ref{eq:3hchi}) and (\ref{eq:invmeas2})
% \begin{align}
%   h_{\mu'}&=2\chi_{\mu'}=\int\log(R^{\#})^2\,d\mu'+\int\log\kappa(Rz)\,d\mu'-\int
%   \log\kappa(z)\,d\mu'   
%   \notag 
%   \\
% &=\int\log J_{\mu'}\,d\mu'.\label{eq:rohlin}
% \end{align}
% This is known as \emph{Rohlin's formula}\index{Rohlin's formula}. It
% is valid in greater generality (see \cite{przyurban}, chapter 1). 
In this section, we show
that $\alpha=\frac{\log N}{\chi} >\frac{2}{\dim_H(\SC)}$, unless $R$
is a Latt\`{e}s map, in which case equality holds. By (\ref{eq:Jensen}),
it is enough to show that
for our given measure $\mu$ the Jacobian is not constant. 

\begin{lemma}
  \label{lem:Jacoinv}
  The Jacobian $J_{\mu}$ of the invariant measure $\mu$ satisfies the
  following.  

  \begin{enumerate}
  \item 
    \label{lem:Jacoinv1}
    It is continuous on $\CDach\setminus(\crit(R)\cup \post(R))$. For
    $\dist(x,\crit(R)\cup\post(R))\geq \epsilon$ 
    \begin{equation*}
      J_\mu(x)\asymp 1,
    \end{equation*}
    where $C(\asymp)=C(\epsilon)$. 
%   \item \label{lem:Jacoinv2}At a postcritical point $p$ with multiplier $\Lambda=\Lambda(p)=R^{\#}(p)$:
%     $$J_{\mu'}(p)=\Lambda^{\frac{2}{\max\deg(p)}}.$$
  \item 
    \label{lem:Jacoinv3}
    In the neighborhood of a point $q\in \crit(R)\cup \post(R)$, with
    $p=R(q)$ 
    \begin{equation*}
      J_{\mu}(z)\asymp
      \abs{z-q}^{-2\left(\frac{1}{\max\deg(q)}-\frac{\deg(q)}{\max\deg(p)}\right)}. 
    \end{equation*}
  \end{enumerate}
\end{lemma}

Note that $\max\deg(p) \geq \deg(q) \max\deg (q)$, so the exponent in
the last expression is nonpositive. 

\begin{proof}
  (\ref{lem:Jacoinv1}) is clear from (\ref{eq:2defJac}) and Lemma
  \ref{lem:kjtok} (\ref{lemkasympcont}), (\ref{lem:kasymp1}). 

  \medskip
  (\ref{lem:Jacoinv3}) Let $z$ be contained in a suitably small
  neighborhood of $q\in \crit(R)\cup\post(R)$. Let $p=R(q)\in
  \post(R)$, $d=\deg_R(q)$. Note that $d$ or $\max\deg(q)$ may be $1$,
  though not both at the same time.  
  We obtain (using Lemma \ref{lem:kjtok})
  \begin{align*}
    & R^{\#}(z) \asymp \abs{z-q}^{d-1}
    \\
    & \kappa(R(z)) \asymp
    \abs{R(z)-p}^{-2\left(1-\frac{1}{\max\deg(p)}\right)}
    \\
    & \phantom{\kappa(R(z))} \asymp \abs{z-q}^{-2d\left(1-\frac{1}{\max\deg(p)}\right)}
    \\
    & \kappa(z) \asymp
    \abs{z-q}^{-2\left(1-\frac{1}{\max\deg(q)}\right)}.
    \intertext{Thus}
    & J_\mu(z) = R^{\#}(z)^2\frac{\kappa(R(z))}{\kappa(z)}\asymp
    \abs{z-q}^{-2\left(\frac{1}{\max\deg(q)}-\frac{\deg(q)}{\max\deg(p)}\right)},
  \end{align*}
  as desired.

%   Let us look at the behavior at a postcritical point $p$ first. By examining the last proof, there is a constant $a\neq 0$ such that
%   \begin{align*}
%     &\kappa(z)\asymp a\dist(z,p)^{-2\left(1-\frac{1}{\max\deg(p)}\right)},
%     \\
%     \intertext{as well as }
%     &\dist(Rz,p)\asymp\Lambda(p)\dist(z,p),
%   \end{align*}
%   where $C(\asymp)\to 1$ as $z\to p$.
%    Plugging this into equation (\ref{eq:Jeq}) yields
%   \begin{equation*}
%     J_{\mu'}(p)=\Lambda^2\Lambda^{-2\left(1-\frac{1}{\max\deg(p)}\right)}=\Lambda^{\frac{2}{\max\deg(p)}}.
%   \end{equation*}
%   Consider now the behavior near a critical point $c$ which gets mapped to $p$. Since
%   \begin{equation*}
%     \dist(z,c)^{\deg(c)}\asymp\dist(Rz,p),
%   \end{equation*}
%   we get by property \ref{lem:sasymppole} of lemma \ref{lem:densitybehavior}
%   \begin{align*}
%     \kappa(Rz)&\asymp\dist(Rz,p)^{-2\left(1-\frac{1}{\max\deg(p)}\right)}
%     \\
%     &\asymp\dist(z,c)^{-2\deg(c)\left(1-\frac{1}{\max\deg(p)}\right)}.
%   \end{align*}
%   Also $\kappa(z)\asymp 1$ and $R^{\#}(z)\asymp \dist(z,c)^{\deg(c)-1}$ and so by formula (\ref{eq:Jeq})
%   \begin{equation*}
%     J_{\mu'}(z)\asymp\dist(z,c)^{-2\left(1-\frac{\deg(c)}{\max\deg(p)}\right)}.
%   \end{equation*}
%   The continuity away from critical points is clear.

\end{proof}

\begin{lemma}
  \label{lem:alpha}
  It holds $\alpha> \frac{2}{\dim_H(\SC)}$ unless $R$ is a Latt\`{e}s
  map, in which case $\alpha = \frac{2}{\dim_H(\SC)}$.
\end{lemma}

\begin{proof}
  The reader should review the terminology from Section
  \ref{sec:orbif-assoc-r}. 
  From (\ref{eq:Jensen}) and Theorem \ref{thm:dS2}
  (\ref{eq:1dimStilde}) it follows that $\alpha = \frac{\log N}{\chi}
  \geq \frac{2\log N}{\log \deg R} = \frac{2}{\dim_H(\mathcal{S})}$. 
  Assume now that
  $\alpha=\frac{2}{\dim_H(\SC)}$, which is equivalent to 
  $J_\mu=\deg R$ almost everywhere. This means in
  particular that the exponent in Lemma \ref{lem:Jacoinv}
  (\ref{lem:Jacoinv3}) is $0$, or $\max\deg(p)=\deg(q)\max\deg(q)$
  (for $R(q)=p$).  
  Thus, for all $x,y\notin \post(R)$ with $R^n(x)=R^n(y)$ we have
  $\deg_{R^n}(x)=\deg_{R^n}(y)$. This implies that the ramification
  function equals the maximal degree, $\nu(p)=\max\deg(p)$ for all
  $p\in \post(R)$. 
  Let $k=\#\post(R)$, and $(\nu(p_1),\dots, \nu(p_k))$ be the
  signature of the 
  orbifold associated to $R$. 

  \medskip
  The argument will be a simple counting argument using the Euler
  characteristic. Fix a $j\geq 0$. 
  Consider the $j$-tiles, $j$-edges, and $j$-vertices;
  they are called faces, edges, vertices in the following for
  simplicity. 
  Let $F,E,V$ be the number of faces, edges, vertices. Since every
  face is incident to $k$ edges, and every 
  edge is incident to two faces we have
  \begin{equation*}
    kF=2E.
  \end{equation*}
  
  Consider now a vertex $q$ that is mapped by $R^j$ to (the first
  postcritical point) $p_1$. Then $q$
  is contained in $2\deg_{R^j}(q)$ faces. 
  Note that each face
  contains exactly one such vertex $q$. Let $n_1$ be the
  number of such vertices $q$. 
  By the above
  $\deg_{R^j}(q)=\nu(p_1)$ if $q\notin \post(R)$, otherwise
  $\deg_{R^j}(q)\leq \nu(p_1)$. Therefore, 
  \begin{equation*}
    2(n_1 -k) \nu(p_1) \leq F \leq 2n_1 \nu(p_1).
  \end{equation*}
  Analog expressions of course hold for all $p_l\in \post(R)$. 
  Thus,
  \begin{equation*}
    \frac{F}{2} \sum_l\frac{1}{\nu(p_l)} 
    \leq
    \sum_l n_l = V 
    \leq 
    \frac{F}{2} \sum_l\frac{1}{\nu(p_l)} + k^2.  
  \end{equation*}
  This yields
  \begin{align*}
    2 & = F-E+V = F - \frac{k}{2} F + \frac{1}{2} F \sum_l
    \frac{1}{\nu(p_l)} + O(1)
    \\
    & = \frac{F}{2}\left(2 - k +
      \sum_l\frac{1}{\nu(p_l)}\right) + O(1)
    \\
    & = \frac{F}{2} \underbrace{\left(2 + \sum_l \left( \frac{1}{\nu(p_l)} -
        1\right)\right)}_{\chi(\OC)} + O(1).
  \end{align*}
  
  Since $F= 2(\deg R)^j$, this can only be satisfied if $\chi(\OC)=0$. This
  means that the orbifold is parabolic, thus $R$ is a Latt\`{e}s map. 
  
\end{proof}

% \begin{remark}
%   In the case 
% \end{remark}

% \bigskip
% So the Jacobian has poles at the critical points unless for every postcritical point $p$ the degree at all critical points $c\in\{R^{-1}(p)\}$ is the same. This only happens if our self similar snowsphere $\SC$ is just the surface of a cube (meaning the generator is just a square subdivided into small squares, without any ``hills'' or ``valleys''). The rational maps corresponding to euclidean subdivisions (subsection \ref{sec:eucl-subd}) also satisfy this. Let $M$ be the number of $1/N$-squares in the $N$-generator. Then $\deg R=M$ as well. If the snowsphere $\SC$ is not the surface of the cube (i,~e.\ $M> N^2$) the Jacobian is not constant by the above, hence by Rohlin's formula (\ref{eq:rohlin}) and Jensen's inequality (\ref{eq:Jensen}) 
% \begin{equation*}
%   h_{\mu'}<\log M.
% \end{equation*}
% Recall that the dimension of the elliptic harmonic measure $\mu_f$
% given by $f\colon \SC\to S^2$ is the same as the dimension of the measure $\mu:=\mu'\circ \tilde{f}$ given by $\tilde{f}\colon\widetilde{\SC}\to\Hob$. Since $h_{\mu}=h_{\mu'}$ we have by the above, Manning's formula (\ref{eq:manning}), and equation (\ref{eq:lyap2})
% \begin{equation*}
%   \dim\mu=\frac{h_{\mu}}{\chi_\mu}<\frac{\log M}{\log N}=\dim_H\SC.
% \end{equation*}
% We have proved the following.
% \begin{theorem}
%   \label{thm:dimell}
%   The dimension of an elliptic harmonic measure on a self similar snowball is strictly less than the Hausdorff dimension of the snowsphere $\SC$, unless $\SC=\partial [0,1]^3$.
% \end{theorem}

\section{Numerical Experiments}
\label{sec:numer-exper}
In principle,
one could use the explicit description of the invariant measure $\mu$
to calculate the entropy or the Lyapunov exponent by equation
(\ref{eq:Rohlinsform}) or (\ref{eq:defLyap}) and thus the dimension of
elliptic harmonic measure. However, we do not know the rate of 
convergence of the densities $\kappa_j\to\kappa$. We use Birkhoff's
ergodic theorem instead to calculate the Lyapunov exponent 
\begin{align}
  \label{eq:lyap1a} 
  \frac{1}{n}&\sum_{j=0}^{n-1}\log R^{\#}(R^jz_0)\to \chi,
  \text{ as } n\to\infty, 
%   \\
%   \frac{1}{n}&\sum_{j=0}^{n-1}\log \abs{R'(R^jz_0)}\to \chi_{\mu'},
%   \text{ as } n\to\infty,\label{eq:lyap1b} 
\end{align}
for ($\mu$ or Lebesgue) almost every $z_0$.
%  \frac{1}{j}\log\abs{(R^j)^{\#}(z_0)}\to \chi_{\mu'}, \text{ for } \mu'\text{-a.~e.\ } z_0.
The dimension of $\mu$ with respect to the metric $\abs{x-y}_{\SC}$ is
then $\dim \mu =\frac{2\chi}{\log N}$. 
% From this one gets the entropy of the elliptic harmonic measure $h_{\mu}=h_{\mu'}=2\chi_{\mu'}$ and its dimension $\dim \mu=\frac{h_{\mu}}{\chi_{\mu}}=\frac{h_{\mu}}{\log N}$.
The ergodic sum is very easy to calculate. There is no way to
determine however that we picked a \emph{generic point} $z_0$. Even if
we did, we do not know how fast the sum converges. So the results in
this section are to be understood as \emph{numerical experiments}
only. Still the values should give a good indication of what to
expect.  

The dimensions thus found are recorded in 
Table \ref{tab:RS}.  
We picked $100$ random
starting values $z_0\in [0,1]^2$, uniformly distributed. The number of
iterations was in each case $n=10000$. 
The computations are done for
our standard example (see Section \ref{sec:orig-with-rati})
$\widehat{R}$. The other examples are from \cite{snowemb}. 
The map $R$
embeds a snowsphere where the generator is a square divided into
$3\times 3$ squares, with a $1/3$-cube put on the middle $1/3$-square. The
maps $R_6,R_7$ embed ``triangular''
snowspheres. More precisely, the generator is a unit equilateral triangle
divided into $4$ equilateral triangles of side-length $1/2$. The
snowsphere represented by $R_6$ has a generator where the middle
triangle is replaced by a tetrahedron. The snowsphere represented by
$R_7$ has a generator where the middle triangle is replaced by a
octahedron. The snowspheres represented by $R,R_6,R_7$ do have self
intersections when embedded in $\R^3$. 

The maps $R_1,R_2,R_3,R_4$ are Latt\`{e}s maps, and serve as a
``control group'' for our computations. Their signatures are
$(2,4,4)$, $(3,3,3)$, $(2,2,2,2)$, and $(2,3,6)$. 
%   from subsection \ref{subsec:varsnow} represented by $R_6$
% and $R_7$. Recall that these are abstract snowspheres, not embedded
% ones. This explains that the dimension of the snowsphere represented
% by $R_7$ is bigger than $3$.   
% The snowsphere represented by $R_8$ was tried as well, but the numeric
% seems to fail in this case (negative dimensions were obtained). For
% comparison the experiments were done for the euclidean subdivisions
% represented by $R_1,R_2,R_3$, and $R_4$. Here the dimension of the
% elliptic harmonic measure is $2$, so these examples serve as a test
% for our method. We record the average dimension from our sample,
% maximal and minimal dimension of the calculated elliptic harmonic
% measure, and the standard deviation. 
The computations were done with
Maple 9, using machine precision. 
\begin{table}[htbp]
  \centering
  \caption{Experimental values of the dimensions by equation (\ref{eq:lyap1a}).}
  \begin{tabular}{llllll}
    \toprule
    rational & dimension of & $\dim \mu$  & $\dim \mu$ & $\dim \mu$
    & standard \\[-3pt] 
    map      & snowsphere   & average     & maximal    & minimal &
    deviation 
    \\
    \midrule
    \\
    [-12pt]
    $\widehat{R}$ & $\frac{\log 29}{\log 5}=2.092\dots$ & $2.030$ &
    $2.034$ & $2.026$ & $0.0017$ 
    \\
    [4pt]
    ${R}$ & $\frac{\log 13}{\log 3}=2.335\dots$ & $2.115$ & $2.125$ &
    $2.105$ & $0.0044$
    \\
    [4pt]
    $R_6$         & $\frac{\log 6}{\log 2}=2.585\dots$ & $2.359$ &
    $2.376$ & $2.342$ & $0.0068$ 
    \\
    [4pt]
    $R_7$         & $\frac{\log 10}{\log 2}=3.322\dots$ & $2.594$ &
    $2.625$ & $2.563$ & $0.012$
    \\[5pt]
    \midrule
    \\[-12pt]
    $R_1$         & $2$        & $1.999$ & $2.001$ & $1.996$ &
    $0.00079$
    \\
    [4pt]
    $R_2$         & $2$        & $1.9998$ & $2.0003$ & $1.9987$ &
    $0.00025$
    \\
    [4pt]
    $R_3$         & $2$        & $1.9999$ & $2.0003$ & $1.9996$ &
    $0.00012$
    \\
    [4pt]
    $R_4$         & $2$        & $1.9997$ & $2.0023$ & $1.9943$ &
    $0.0012$
    \\
    \bottomrule
  \end{tabular}

  \label{tab:RS}
\end{table}

\section{Open Problems}
\label{sec:open-problems}

We conclude with some questions. 
\begin{problem}
  Is it true that
  \begin{equation*}
    \dim \mu>2 ?
  \end{equation*}
  Here, $\mu$ is an elliptic harmonic measure on a self similar snowsphere.
\end{problem}
The answer is generally expected to be yes. For corresponding results
on the dimension of harmonic measure see \cite{Wolff}. 
There is some hope that the description through the
rational map $R$ might help answering this question. 
\begin{problem}
  Is there a nontrivial upper bound for the dimension of an elliptic
  harmonic measure? More precisely, is there an $\epsilon>0$, such that  
  \begin{equation*}
    \dim \mu\leq 3-\epsilon,
  \end{equation*}
  for an elliptic harmonic measure $\mu$ of any snowsphere
  $\SC \subset\R^3$?
\end{problem}
The corresponding statement for harmonic measure is true by
\cite{Bour}, see also \cite{Mak} for the two-dimensional case. 

\begin{problem}
  Is there a \emph{geometric way} to ``see'' the exponent $\alpha$/the
  dimension of the elliptic harmonic measure? 
  This would characterize the ``quasisymmetric size'' of a quasisphere. 
  In general, one expects
  to have different upper and lower dimensions ($\dim^*\mu = \dim \mu$ as
  before, $\dim_* \mu = \inf\{\dim_H K \mid \mu(K)> 0$\}).
\end{problem}

\begin{problem}
  Given the elliptic harmonic measure on a snowsphere $\SC\subset
  \R^3$, is there a \emph{relation to other 
  natural measures} on $\SC$? In particular, is there a relation to the
  harmonic measure or the $3$-harmonic measure on $\SC$?
\end{problem}
% \begin{table}[htbp]
%   \centering
%   \caption{Experimental values of the dimensions by equation (\ref{eq:lyap1b})}
%   \begin{tabular}{llllll}
%     \toprule
%     rational & dimension of & $\dim \mu$  & $\dim \mu$ & $\dim \mu$   & standard \\[-6pt]
%     map      & snowsphere   & average     & maximal    & minimal & deviation \\
%     \midrule
%     $\widehat{R}$ & $\frac{\log 13}{\log 3}=2.335\dots$    & $2.115$ & $2.126$ & $2.105$ & $0.0044$\\
%     $R_6$         & $\frac{\log 6}{\log 2}=2.585\dots$    & $2.359$ & $2.376$ & $2.344$ & $0.0067$\\
%     $R_7$         & $\frac{\log 10}{\log 2}=3.322\dots$    & $2.597$ & $2.625$ & $2.572$ & $0.010$\\[5pt]
%     \midrule
%     $R_1$         & 2        & $2.000$ & $2.004$ & $1.996$ & $0.0014$\\
%     $R_2$         & 2        & $2.000$ & $2.004$ & $1.998$ & $0.00070$\\
%     $R_3$         & 2        & $2.000$ & $2.002$ & $2.000$ & $0.00049$\\
%     $R_4$         & 2        & $2.000$ & $2.003$ & $1.992$ & $0.0015$\\
%     \bottomrule
%   \end{tabular}
%   \label{tab:RE}
%\end{table}

% ==========   Bibliography
%
\bibliographystyle{amsalpha}
\bibliography{litlist}

\end{document}